\documentclass[final,leqno,onefignum]{siamltex704}

\usepackage{mathrsfs}
\usepackage{amsmath,amssymb}
\usepackage{enumitem}
\usepackage{amsfonts}
\usepackage{latexsym}
\usepackage{pictex}
\usepackage{color}
\usepackage{graphicx}

\newtheorem{remark}{Remark}[section]
\newtheorem{assumption}{Assumption}[section]

\def\lb{\label}

\def\bf{\textbf}

\def\g{\gamma}
\def\mr{{\mathbb R}}
\def\mn{{\mathbb N}}
\def\R{{\Bbb R}}
\def\Z{{\Bbb Z}}
\def\p{\partial}
\def\sign{\mathop{\rm sign}\nolimits}
\def\ha{{\frac{1}{2}}}

\newcommand{\ii}{\mathrm{i}}
\newcommand{\pdpd}[2]{\frac{\partial #1}{\partial #2}}
\newcommand{\vphi}{\varphi}

\begin{document}

\date{\today}

\title{Semiclassical inverse spectral problem for elastic Love waves
  in isotropic media}

\author{Maarten V. de Hoop \thanks{Simons Chair in Computational and
    Applied Mathematics and Earth Science, Rice University, Houston,
    TX 77005, USA (mdehoop@rice.edu)}
\and Alexei Iantchenko \thanks{Department of Materials Science and
  Applied Mathematics, Faculty of Technology and Society, Malm\"{o}
  University, SE-205 06 Malm\"{o}, Sweden (ai@mau.se)}
\and Robert D. van der Hilst \thanks{Department of Earth, Atmospheric
  and Planetary Sciences, Massachusetts Institute of Technology,
  Cambridge, MA 02139, USA (hilst@mit.edu)}
\and Jian Zhai \thanks{Institute for Advanced Study, The Hong Kong University of Science and Technology, Kowloon, Hong Kong, China (jian.zhai@outlook.com)}}

\maketitle

\pagestyle{myheadings}
\thispagestyle{plain}
\markboth{DE HOOP, IANTCHENKO, VAN DER HILST and ZHAI}{Semiclassical
  inverse spectral problem for Love waves}

\begin{abstract}
We analyze the inverse spectral problem on the half line associated
with elastic surface waves. Here, we focus on Love waves. Under
certain generic conditions, we establish uniqueness and present a
reconstruction scheme for the \textit{S}- wavespeed with multiple
wells from the semiclassical spectrum of these waves.
\end{abstract}

\section{Introduction}
\label{intro}

We analyze the inverse spectral problem on the half line associated
with elastic surface waves. Here, we focus on Love waves. In a
follow-up paper we present the corresponding inverse problem for
Rayleigh waves. Surface waves have played a key role in revealing
Earth's structure from the shallow near-surface to several hundred
kilometers deep into the mantle, depending on the frequencies and data
acquisition configurations considered.

\subsection{Seismology}

The inverse spectral problem for surface waves fits in the
seismological framework of surface-wave tomography. Surface-wave
tomography has a long history. Since pioneering work on inference from
the dispersion of surface waves half a century ago \cite{Haskell,
  Press, BD, Knopoff, Toksoz, Schwab, Dziewonski, Woodhouse, Nolet1},
surface wave tomography based on dispersion of waveforms from
earthquake data has played an important role in studies of the
structure of the Earth's crust and upper mantle on both regional and
global scales \cite{Lerner, Woodhouse2, Nataf, Nolet2, MT, Trampert,
  GJG, Ritzwoller-2001, Simons, BE, Ritsema-2004, Lebedev-2008, Yao2}.

In order to avoid the effects of scattering due to complex crustal
structure, these studies focused on the analysis, measurement, and
inversion of surface wave dispersion at relatively low frequencies
(that is, $4-20$ mHz, or periods between $50$ to $250$ s) at which the
fundamental modes sense mantle structure to $200-300$ km depth and
higher modes reach across the upper mantle and transition zone to some
$660$ km depth. Most methods assume some form of (WKB) asymptotic and
path-average approximation \cite{DT} in line with our semiclassical
point of view.

More than a decade ago, Campillo and his collaborators discovered that
cross correlation of ambient noise yields Green's function for surface
waves \cite{Campillo1,Campillo2,Campillo3}. This enabled the
possibility to extend the applicability of surface-wave tomography not
only to any area where seismic sensors can be placed, but also to
short-path measurements and frequencies at which the data are most
sensitive to shallow depths.  Crustal studies based on ambient noise
tomography are typically conducted in the period band of $5–40$ s, but
shorter period surface waves ($\sim 1$ s, using station spacing of
$\sim 20$ km or less) have been used to investigate shallow crustal or
even near surface shear-wave speed variations \cite{Sabra, Yao1, Yao2,
  Yang, Lin, Huang}.

\subsection{Semiclassical analysis perspective}

In a separate contribution \cite{dHINZ}, we presented the
semiclassical analysis of surface waves. Such an analysis leads to a
geometric-spectral description of the propagation of these waves
\cite{Babich, Woodhouse}. This semiclassical analysis is built on the
work of Colin de Verdi\`{e}re \cite{CDV2, CDV3}. \textit{The main
  contribution of this paper is the construction of the
  Bohr-Sommerfeld quantization for Love waves.} Colin de Verdi\`{e}re
also considered the inverse spectral problem of scalar surface waves
allowing wavespeed profiles that contain a well \cite{CdV2011}. His
result does not account for the Neumann boundary condition at the
surface, although a reflection principle could be invoked, but his
methodology directly applies once the Bohr-Sommerfeld quantization is
obtained. The reflection principle does not apply to general elastic
surface waves and the remedy is presented in this paper. In the
process, we show that with the Neumann boundary condition at the
surface, in fact, ambiguities arising in the recovery of the
\textit{S}-wave speed on the line (that is, without this boundary
condition) can be resolved.

We study the elastic wave equation in $X = \mr^2 \times
(-\infty,0]$. In coordinates,
$$
   (x,z) ,\quad x = (x_1,x_2) \in \mr^2 ,\
                         z \in \R^{-} = (-\infty,0],
$$
we consider solutions, $u = (u_1,u_2,u_3)$, satisfying the Neumann
boundary condition at $\partial X = \{z=0\}$, to the system
\begin{equation}\label{elaswaeq}
\begin{split}
   \partial^2_t u_i + M_{il} u_l &= 0 ,\\
   u(t=0,x,z) &= 0,~~\partial_tu(t=0,x,z)=h(x,z) ,\\
   \frac{c_{i3kl}}{\rho}\partial_k u_l(t,x,z=0) &= 0 ,
\end{split}
\end{equation}
where
\begin{multline*}
M_{il} =
-\frac{\partial}{\partial z}\frac{c_{i33l}(x,z)}{\rho(x,z)}
\frac{\partial}{\partial z}     - \sum_{j,k=1}^{2}
\frac{c_{ijkl}(x,z)}{\rho(x,z)} \frac{\partial}{\partial x_j}
\frac{\partial}{\partial x_k}
- \sum_{j=1}^{2}
\frac{\partial}{\partial x_j}\frac{c_{ij3l}(x,z)}{\rho(x,z)}
\frac{\partial}{\partial z}
\\
- \sum_{k=1}^{2}
\frac{c_{i3kl}(x,z)}{\rho(x,z)} \frac{\partial}{\partial z}
\frac{\partial}{\partial x_k}
-\sum_{k=1}^{2}
\left( \frac{\partial}{\partial z} \frac{c_{i3kl}(x,z)}{\rho(x,z)} \right)
\frac{\partial}{\partial x_k}
- \sum_{j,k=1}^{2}
\left( \frac{\partial}{\partial x_j} \frac{c_{ijkl}(x,z)}{\rho(x,z)} \right)
\frac{\partial}{\partial x_k}.
\end{multline*}
Here, the stiffness tensor, $c_{ijkl}$, and density, $\rho$, are
smooth and obey the following scaling: Introducing $Z =
\frac{z}{\epsilon}$,
$$
   \frac{c_{ijkl}}{\rho}(x,z)
     = C_{ijkl}\left(x,\frac{z}{\epsilon}\right) ,
                            ~~\epsilon \in (0,\epsilon_0];
$$
$$
   C_{ijkl}(x,Z)=C_{ijkl}(x,Z_I)=C_{ijkl}^I(x),\quad
                            Z \leq Z_I<0.
$$
As discussed in \cite{dHINZ}, surface waves travel along the surface
$z = 0$.

\medskip\medskip

\noindent
The remainder of the paper is organized as follows. In Section
\ref{semi}, we give the formulation of the inverse problems as an
inverse spectral problem on the half line. In Section
\ref{S-decreasing}, we treat the simple case of recovery of a
monotonic profile of wave speed. In Section \ref{BS rule}, we discuss
the relevant Bohr-Sommerfeld quantization, which is the corner stone
in the study of the inverse spectral problem. In Section
\ref{inverse}, we give the reconstruction scheme under generic
assumptions.

\section{Semiclassical description of Love waves}\label{semi}

\subsection{Surface wave equation, trace and the data}

For the convenience of the readers, we briefly summarize the
semiclassical description of elastic surface waves. The leading-order
Weyl symbol associated with $M_{il}$ above is given by
\begin{multline}\label{H0}
H_{0,il}(x,\xi) =
-\frac{\partial}{\partial Z}C_{i33l}(x,Z)
\frac{\partial}{\partial Z}
\\
- \ii\sum_{j=1}^{2}
C_{ij3l}(x,Z) \xi_j
\frac{\partial}{\partial Z}
- \ii\sum_{k=1}^{2}
C_{i3kl}(x,Z) \frac{\partial }{\partial Z}
\xi_k
- \ii \sum_{k=1}^{2}
\left( \frac{\partial}{\partial Z} C_{i3kl}(x,Z) \right)
\xi_k 
\\
+ \sum_{j,k=1}^{2}
C_{ijkl}(x,Z) \xi_j \xi_k.
\hspace*{4.0cm}
\end{multline}
We view $H_0(x,\xi)$ as ordinary differential
operators in $Z$, with domain
$$
   \mathcal{D} = \left\{ v \in H^2(\mathbb{R}^-)\ \bigg|\
      \sum_{l=1}^3\left(C_{i33l}(x,0)
    \frac{\partial v_l}{\partial Z}(0)
       + \ii \sum_{k=1}^2C_{i3kl}\xi_kv_l(0)\right) = 0 \right\} .
$$
For an isotropic medium,
$$
   C_{ijkl} = \hat{\lambda} \delta_{ij} \delta_{kl}
       + \hat{\mu} (\delta_{ik} \delta_{jl} + \delta_{il} \delta_{jk}),
$$
where $\hat{\lambda}=\frac{\lambda}{\rho}$ and
$\hat{\mu}=\frac{\mu}{\rho}$.  The \textit{S}-wavespeed, $C_S$, is
then $C_S=\sqrt{\hat{\mu}}$. The decoupling of Love and Rayleigh waves
is observed in practice, and explained in \cite{dHINZ}. We denote
$$
   P(\xi) = \left(\begin{array}{ccc}
        |\xi|^{-1}\xi_2& |\xi|^{-1}\xi_1 & 0 \\
        -|\xi|^{-1}\xi_1 &|\xi|^{-1}\xi_2 & 0 \\
        0 & 0 & 1 \end{array}\right) .
$$
Then 
$$
P(\xi)^{-1}H_0(x,\xi)P(\xi)=\left(\begin{array}{cc}
H_0^L(x,\xi) &\\
&H_0^R(x,\xi)
\end{array}\right),
$$
where
\begin{equation}\label{love}
   H_0^L(x,\xi)\vphi_1 = -\pdpd{}{Z} \hat{\mu} \pdpd{\vphi_1}{Z}
             + \hat{\mu} \, |\xi|^2 \vphi_1
\end{equation}
supplemented with boundary condition
\begin{equation*}\label{love_b}
   \pdpd{\vphi_1}{Z}(0) = 0 ,
\end{equation*}
for Love waves. We will consider only the Love waves in this paper.

We assume that $\Lambda_\alpha(x,\xi)$ is an eigenvalue of
$H_0(x,\xi)$ with eigenfunction $\Phi_{\alpha,0}(Z,x,\xi)$. By
\cite[Theorem 2.1]{dHINZ}, we have
\begin{equation}\label{eq:comm_sys}
   H_0^L \circ \Phi_{\alpha,0}\, 
   = \Phi_{\alpha,0}
     \circ \Lambda_\alpha+\mathcal{O}(\epsilon).
\end{equation}
We define
\begin{equation}\label{defchi}
   J_{\alpha,\epsilon}(Z,x,\xi) =
      \frac{1}{\sqrt{\epsilon}}\Phi_{\alpha,0}(Z,x,\xi) .
\end{equation}
Microlocally (in $x$), we can construct approximate
solutions of the system \eqref{elaswaeq} with initial values
$$
   h(x,\epsilon Z)=\sum_{\alpha=1}^{\mathfrak{M}}
   J_{\alpha,\epsilon}(Z,x,\epsilon D_x) W_{\alpha,\epsilon}(x,Z) ,
$$
representing surface waves. We assume
that all eigenvalues $\Lambda_1 < \cdots < \Lambda_\alpha < \cdots <
\Lambda_{\mathfrak{M}}$ are eigenvalues of the operator given in
(\ref{love}). We let $W_{\alpha,\epsilon} $ solve the initial value
problems (up to leading order)
\begin{eqnarray}
   [\epsilon^2 \p_t^2
       + \Lambda_{\alpha}(x,D_x)]
	 W_{\alpha,\epsilon}(t,x,Z) &=& 0 ,
\label{wqeaj_sys}\\
   W_{\alpha,\epsilon}(0,x,Z) &=& 0 ,\quad
   \p_t W_{\alpha,\epsilon}(0,x,Z) = J_{\alpha,\epsilon}W_{\alpha}(x,Z) ,
\label{eq:initrange}
\end{eqnarray}
$\alpha = 1,\ldots,\mathfrak{M}$. We let
$\mathcal{G}_0(Z,x,t,Z',\xi;\epsilon)$ denote the approximate Green's
function (microlocalized in $x$), up to leading order, for Love
waves. We may write \cite{dHINZ}
\begin{multline}
\label{G0}
\mathcal{G}_0(Z,x,t,Z',\xi;\epsilon)
= \sum_{\alpha=1}^{\mathfrak{M}} J_{\alpha,\epsilon}(Z,x,\xi)
\\
\left(\frac{\mathrm{i}}{2}\mathcal{G}_{\alpha,+,0}(x,t,\xi,\epsilon)-\frac{\mathrm{i}}{2}\mathcal{G}_{\alpha,-,0}(x,t,\xi,\epsilon)\right)\Lambda^{-1/2}_\alpha(x,\xi)J_{\alpha,\epsilon}(Z',x,\xi) ,
\end{multline}
where $\mathcal{G}_{\alpha,\pm,0}$ are Green's functions for half wave
equations associated with (\ref{wqeaj_sys})-(\ref{eq:initrange}).
We have the trace
\[\lb{trace_from_normal_modes}
   \int_{\mathbb{R}^-}
   \widehat{\epsilon
        \partial_t \mathcal{G}_0}(Z,x,\omega,Z,\xi;\epsilon)
        \mathrm{d} \epsilon Z
   = \sum_{\alpha=1}^{\mathfrak{M}}
     \delta(\omega^2 - \Lambda_\alpha(x,\xi))
        \Lambda^{1/2}_\alpha(x,\xi) + {\mathcal O}(\epsilon^{-1})
\]
from which we can extract the eigenvalues
$\Lambda_\alpha,\,\alpha=1,2,\cdots,\mathfrak{M}$ as functions of
$\xi$. We use these to recover the profile of $C_S^2$.

In practice, these eigenvalues are obtained from surface-wave
tomography and to ensure that all eigenvalues are observed,
measurements of surface-waveforms should be taken in boreholes. Most
seismic observations are made at or near Earth’s surface, but modern
networks increasingly include borehole sensors indeed. For example,
Hi-net seismographic network in
Japan\footnote{http://www.hinet.bosai.go.jp} includes more than $750$
sensors located in $> 100$ m deep boreholes and permanent sites of
USArray\footnote{http://www.usarray.org} include sensors placed
around $100$ m depth.

\subsection{Semiclassical spectrum}

From here on, we only consider the operator $H_0^L(x,\xi)$ for Love
waves. We suppress the dependence on $x$, and introduce $h =
|\xi|^{-1}$ as another semiclassical parameter. Within this setting,
we also change the notation from $\frac{\partial}{\partial Z}$ to
$\frac{\mathrm{d}}{\mathrm{d}Z}$. We arrive at the operator
\begin{equation*}
   L_h = -h^2 \frac{\mathrm{d}}{\mathrm{d}Z}
      \left(\hat{\mu}(Z) \frac{\mathrm{d}}{\mathrm{d}Z}\right)
            + \hat{\mu}(Z)
\end{equation*}
with Neumann boundary condition at $Z = 0$. The assumption on the
stiffness tensor gives us the following assumption on $\hat{\mu}$:

\medskip

\begin{assumption} \label{assu1}
The (unknown) function $\hat{\mu}$ satisfies $\hat{\mu}(Z) =
\hat{\mu}(Z_I)$ for all $Z \leq Z_I$ and
\begin{equation*}
   0 < \hat{\mu}(0) = E_0
       = \inf_{Z \le 0} \hat{\mu}(Z)
   < \hat{\mu}_I = \sup_{Z \leq 0} \hat{\mu}(Z) = \hat{\mu}(Z_I) .
\end{equation*}
\end{assumption}

\medskip

\noindent
The assumption that $\hat{\mu}$ attains its mininum at the boundary,
and its maximum in some deep zone, is realistic in practice.

We first observe that the spectrum of $L_h$ is divided in two parts,
$$
   \sigma(L_h) = \sigma_{pp}(L_h) \cup \sigma_{ac}(L_h) ,
$$
where the point spectrum $\sigma_{pp}(L_h)$ consists of a finite
number of eigenvalues in $(E_0,\hat{\mu}_I)$ and the continuous
spectrum $\sigma_{ac}(L_h) = [\hat{\mu}_I,\infty)$. We write
  $\lambda_\alpha = h^2 \Lambda_\alpha$. Since this is a
  one-dimensional problem, the eigenvalues are simple and satisfy
\begin{gather*}
   E_0 < \lambda_1(h) < \lambda_2(h) < \ldots
            < \lambda_{\mathfrak{M}}(h) < \hat{\mu}_I ; 
\end{gather*}
the number of eigenvalues, $\mathfrak{M}$ increases as $h$
decreases.

We will study how to reconstruct the profile $\hat{\mu}$ using only
the asymptotic behavior of $\lambda_\alpha(h)$ in $h$. To this end, we
introduce the semiclassical spectrum as in \cite{CdV2011}

\medskip

\begin{definition}
For given $E$ with $E_0 < E \leq \hat{\mu}_I$ and positive real number
$N$, a sequence $\mu_\alpha(h)$, $\alpha=1,2,\ldots$ is a
semiclassical spectrum of $L_h$ mod $o(h^N)$ in $]-\infty,E[$ if, for
  all $\lambda_\alpha(h) < E$,
$$
   \lambda_\alpha(h) = \mu_\alpha(h) + o(h^N)
$$
uniformly on every compact subset $K$ of $]-\infty,E[\,$.
\end{definition}

\section{Reconstruction of a monotonic profile} \lb{S-decreasing}

In this section, we give a reconstruction scheme for the simple
situation where the profile $\hat{\mu}$ is monotonic. First it is well
known that

\medskip

\begin{lemma} \label{first}
The first eigenvalue of $L_h$ satisfies $\lim_{h \rightarrow
  0} \lambda_1(h) = E_0$.
\end{lemma}

\medskip

Similar to Theorem 3 in \cite{CDV2}, we have

\medskip

\begin{theorem} \label{CdV-decr}
Assume that $\hat{\mu}$ is decreasing in $[Z_I,0]$. Then the
asymptotics of the discrete spectra $\lambda_j(h), $ $1 \leq j \leq
\mathfrak{M}_j$ as $h \to 0$ determine the function $\hat{\mu}$.
\end{theorem}

\medskip

\noindent
Before giving the proof, we recall the Abel transform and its
inverse. We introduce
\[
   \mathcal{A} g(E)
           = \int_{E_0}^E \sqrt{E - u} g(u) \, \mathrm{d}u .
\]
Then
\[
   \frac{\mathrm{d}}{\mathrm{d}E} \mathcal{A} g(E)
              = \frac{1}{2} T g(E) ,
\quad
   T g(E) = \int_{E_0}^E \frac{g(u)}{\sqrt{E - u}}
                       \, \mathrm{d}u ,
\]
where $T g$ signifies the Abel transform of $g$. By the inversion
formula for the Abel transform,
\[
   \frac{\mathrm{d}}{\mathrm{d}E} T^2 g(E) = \pi g(E) ,
\]
we get
\begin{equation} \label{eq:Atinv}
   \left(\frac{4}{\pi} \frac{\mathrm{d}^2}{\mathrm{d}E^2}
       \mathcal{A} \frac{\mathrm{d}}{\mathrm{d}E}
       \mathcal{A}\right) g(E) = g(E) .
\end{equation}

\medskip

\begin{proof}
First, we note that $E_0 = \hat{\mu}(0)$ is determined by the first
  semiclassical eigenvalue $\lambda_1(h)$ by Lemma~\ref{first}. Then,
  we invoke Weyl's law. For $E < \hat{\mu}_I$, let $N(h,E) = \#\{
  \lambda_j(h) \leq E \}$, where $\lambda_j(h)$ is an eigenvalue for
  $L_h$. Then \cite{dHINZ}
\begin{equation} \label{countingf}
   N(h,E) = \frac{1}{2\pi h} \left[
   {\rm Area}(\{(Z,\zeta)\ :\ \hat{\mu}(Z) (1 + \zeta^2) \leq E\})
            + o(1)\right] .  
\end{equation}
Thus, from the leading order behavior (in $h$) of $\lambda_j(h)$ we
can recover
\[
   {\rm Area}(\{(Z,\zeta)\ :\ \hat{\mu}(Z) (1 + \zeta^2) \leq E\})
   = 2 \tilde{S}^1_0(E) ,\quad
   \tilde{S}^1_0(E) = \int_{f(E)}^0
            \sqrt{\frac{E - \hat{\mu}}{\hat{\mu}}} \, \mathrm{d}Z ,
\]
with $\hat{\mu}(f(E)) = E$. We change variable of integration, $Z =
f(u)$, with
\begin{equation} \label{sC}
   \left.\frac{\mathrm{d}}{\mathrm{d}Z} \hat{\mu}(Z)
                         \right|_{Z = f(u)} = \frac{1}{f'(u)}
\end{equation}
and get
\[
   \tilde{S}^1_0(E) = \mathcal{A} g(E) ,\quad
                                g(u) = \frac{f'(u)}{\sqrt{u}} .
\]
Applying (\ref{eq:Atinv}) above, we recover $g$, that is,
\[
   f'(E) = \left(\frac{4}{\pi} \sqrt{E}
        \frac{\mathrm{d}^2}{\mathrm{d}E^2} \mathcal{A}
        \frac{\mathrm{d}}{\mathrm{d}E} \right)
                          \tilde{S}^1_0(E) ,\quad
   E_0 < E < \hat{\mu}_I .
\]
Then
\[
   f(E) = \int_{E_0}^E f'(u) \, \mathrm{d}u ,
\]
using that $f(E_0) = 0$ and knowledge
  of $E_0 = \hat{\mu}(0)$ from the first eigenvalue (Lemma
  \ref{first}), from which we recover $\hat{\mu}$ by the inverse
function theorem.
\end{proof}

\section{Bohr-Sommerfeld quantization}
\label{BS rule}

The Bohr-Sommerfeld rules give a quantization for the semiclassical
spectrum \cite{CdV2005}. We will derive these rules making use of the
WKB-Maslov Ansatz for the eigenfunctions. We obtain an alternative
proof to the one given in \cite{CDV3, CdV2011}, which enables to
explicitly incorporate Neumann boundary conditions at the surface. It
opens the way for studying inverse problems also for Rayleigh waves;
these will be investigated in a follow-up paper.

We construct WKB solutions of the form
\begin{equation} \label{WKB}
   u_h(Z) = C \exp\left[\frac{1}{h}
              \sum_{j=0}^\infty h^j \mathcal{S}_j(Z) \right]
\end{equation}
that satisfy
\begin{equation}\label{eq10}
   -h^2 \hat{\mu}(Z) u_h''(Z) - h^2 \hat{\mu}'(Z) u_h'(Z)
        + \hat{\mu}(Z) u_h(Z) = E u_h(Z) .
\end{equation}
We will follow various calculations from \cite{Bender} in the
following analysis.

\subsection{Half well}

We consider the eigenvalue problem on the half line $\mathbb{R}^-$,
with Neumann boundary condition at $Z = 0$. We fix a real number $E$
and assume that there exists a unique $Z_E$ such that $\hat{\mu}(Z_E) =
E$. For exposition of the construction, we change the variable $Z
\rightarrow Z_E - Z$ such that $\hat{\mu}(0) = E$ and $Z_E$ is the
boundary point. Furthermore, we assume that $\hat{\mu}(Z) - E > 0$ for
$Z > 0$ and $\hat{\mu}(Z) - E < 0$ for $Z_E < Z < 0$. The original
domain $]-\infty,0]$ changes to $[Z_E,\infty[$. We divide the domain
    $[Z_E,\infty[$ into three regions: region $\mathrm{I}$ $(Z > 0)$,
region $\mathrm{II}$ ($|Z|$ is small) and region $\mathrm{III}$ $(Z_E
\leq Z < 0)$. We will construct WKB solutions in each region and glue
them together.

First, we construct the WKB solution, $u_\mathrm{I}(Z)$, in
\textbf{region I}. We substitute solutions of the form (\ref{WKB}),
collect terms of equal orders in $h$, and arrive at an infinite family
of equations which may be solved recursively. The $\mathcal{O}(h^0)$
terms give the eikonal equation for $\mathcal{S}_0$,
$$
   \hat{\mu}(Z) (1 - (\mathcal{S}_0'(Z))^2) = E .
$$
We select the solution
\begin{equation}
   \mathcal{S}_0(Z) = -\int_0^Z
      \sqrt{\frac{\hat{\mu} - E}{\hat{\mu}}} \mathrm{d}Z' .
\end{equation}
Then the $\mathcal{O}(h)$ term yields
$$
   \hat{\mu} \mathcal{S}_0''
           + 2 \hat{\mu} \mathcal{S}_0' \mathcal{S}_1'
                 + \hat{\mu}' \mathcal{S}_0' = 0 ,
$$
which implies that
$$
   \mathcal{S}_1' = -\frac{1}{2} (\log(\hat{\mu} \mathcal{S}_0'))'
     = -\frac{1}{4} \left(\log[\hat{\mu} (\hat{\mu} - E)]\right)' ;
$$
we select the solution
\begin{equation}
   \mathcal{S}_1 = -\frac{1}{4} \log[\hat{\mu} (\hat{\mu}  -E)] .
\end{equation}
The lower order terms give us a sequence of equations,
$$
   2 \hat{\mu} \mathcal{S}_0' \mathcal{S}_j'
       + (\hat{\mu} \mathcal{S}_{j-1}')'
     + \hat{\mu} \sum_{k=1}^{j-1} \mathcal{S}_{j-k}' \mathcal{S}_k'
   = 0 ,\quad j \geq 2 .
$$
We write down the explicit form of $\mathcal{S}_2$ for later use
\begin{equation} \label{eq:calS2}
   \mathcal{S}_2(\delta,Z) = \int_\delta^Z
     \left[\frac{(E \hat{\mu}' - 2 \hat{\mu} \hat{\mu}')^2}{
                 32 \hat{\mu}^{3/2} (\hat{\mu} - E)^{5/2}}
   + \frac{-E^2 \hat{\mu}'' + 3 E \hat{\mu} \hat{\mu}''
          - 2 \hat{\mu}^2 \hat{\mu}'' + E (\hat{\mu}')^2}{
   8 (\hat{\mu} - E)^{5/2} \hat{\mu}^{1/2}}\right] \mathrm{d}Z' ,
\end{equation}
up to a constant difference; here $\delta$ is any small fixed positive
constant. Upon integrating by parts, we obtain
\begin{equation}
\begin{split}
   \mathcal{S}_2(\delta,Z)
   = -\frac{(3 E + 2 \hat{\mu}) \hat{\mu}'}{48 \hat{\mu}^{1/2}
                     (\hat{\mu} - E)^{3/2}}
         -&\frac{\hat{\mu}'}{24 (\hat{\mu} - E)^{1/2} \hat{\mu}^{1/2}}
\\
   &+ \int_\delta^Z \left[-\frac{(\hat{\mu}')^2}{24 \hat{\mu}^{3/2}
           (\hat{\mu} - E)^{1/2}}
     + \frac{(7 E - 8 \hat{\mu}) \hat{\mu}''}{
        48 \hat{\mu}^{1/2}(\hat{\mu} - E)^{3/2}}\right] \mathrm{d}Z' .
\end{split}
\end{equation}

Next, we consider \textbf{region II} containing the turning
point. When $|Z|$ is small, we expand
$$
   \hat{\mu}(Z) - E = a_1 Z + a_2 Z^2 + a_3 Z^3 + \cdots .
$$
Here, $a_1 > 0$. We write $u_{\mathrm{II}}(Z) = \hat{\mu}^{-1/2}(Z) \,
v_{\mathrm{II}}(Z) = (E + a_1 Z + a_2 Z^2 + a_3 Z^3 + \cdots)^{-1/2}
v_{\mathrm{II}}(Z)$. Then we obtain
$$
   -h^2 \frac{\mathrm{d}}{\mathrm{d}Z}
   \left(\hat{\mu}
               \frac{\mathrm{d} u_{\mathrm{II}}}{\mathrm{d}Z}\right)
   = -h^2 \frac{\mathrm{d}}{\mathrm{d}Z}
      \left(\hat{\mu} \frac{\mathrm{d}}{
          \mathrm{d}Z} \hat{\mu}^{-1/2}(Z) v_{\mathrm{II}}(Z)\right)
   = -h^2 \hat{\mu}^{1/2} v_{\mathrm{II}}''
                       + h^2 (\hat{\mu}^{1/2})'' v_{\mathrm{II}} .
$$
Thus, by (\ref{eq10}), we have the following equation for
$v_{\mathrm{II}}$,
\begin{equation} \label{eq11}
   h^2 v_{\mathrm{II}}'' = \left( 1 - E \hat{\mu}^{-1}
     + h^2 \frac{(\hat{\mu}^{1/2})''}{\hat{\mu}^{1/2}} \right)
           \, v_{\mathrm{II}} .
\end{equation}
We further employ the simple asymptotic expansion,
$$
   1 - E \hat{\mu}^{-1}(Z) = b_1 Z + b_2 Z^2 + \cdots ,
$$
where $b_1 = \frac{a_1}{E}$ and $b_2 = \frac{a_2 E -
  a_1^2}{E^2}$. Temporarily, we introduce the scaling $Z = h^{2/3}
b_1^{-1/3} Y$. With abuse of notation for $v_{\mathrm{II}}$,
(\ref{eq11}) gives 
$$
   -h^2 h^{-4/3} b_1^{2/3}
      \frac{\mathrm{d}^2 v_{\mathrm{II}}}{\mathrm{d}Y^2}
   = (b_1 h^{2/3} b_1^{-1/3} Y + b_2 h^{4/3} b_1^{-2/3} Y^2
      + \cdots) \, v_{\mathrm{II}} ,
$$
which can be simplified to
\begin{equation} \label{eq12}
   \frac{\mathrm{d}^2 v_{\mathrm{II}}}{\mathrm{d}Y^2}
        \sim (Y + h^{2/3} b_1^{-4/3} b_2 Y^2) \, v_{\mathrm{II}} ,
\end{equation}
keeping the second-order approximation. We then seek an approximate
solution of the form
$$
   v_{\mathrm{II}}(Y) \sim (1 + \alpha_1 h^{2/3} Y) \,
              \mathrm{Ai}(Y + \beta_1 h^{2/3} Y^2) ,
$$
where $\mathrm{Ai}$ is the Airy function and $\alpha_1$ and $\beta_1$
are constants to be determined. By tedious calculations, we find that
\begin{multline*}
   \frac{\mathrm{d}^2 v_{\mathrm{II}}}{\mathrm{d}Y^2}
   \sim D \Big[
        \alpha_1 h^{2/3} \mathrm{Ai}'(Y + \beta_1 h^{2/3} Y^2)
        + \alpha_1 h^{2/3} (1 + 2 \beta_1 h^{2/3} Y)
          \, \mathrm{Ai}'(Y + \beta_1 h^{2/3} Y^2)
\\
   + (1 + \alpha_1 h^{2/3} Y) (1 + 2 \beta_1 h^{2/3} Y)^2
          \mathrm{Ai}''(Y + \beta_1 h^{2/3} Y^2)
\\
   + (1 + \alpha_1 h^{2/3} Y) \,
          2 \beta_1 h^{2/3} \mathrm{Ai}'(Y + \beta_1 h^{2/3} Y^2)
     \Big] .
\end{multline*}
Comparing this equation with differential equation (\ref{eq12}), and
using that
$$
   \mathrm{Ai}''(Y + \beta_1 h^{2/3} Y^2) = (Y
      + \beta_1 h^{2/3} Y^2) \, \mathrm{Ai}(Y + \beta_1 h^{2/3} Y^2) ,
$$
we must have
$$
   \alpha_1 + \beta_1 = 0 ,
$$
and
$$
   5 \beta_1 = b^{-4/3}_1 b_2 .
$$
Hence, undoing the scaling and returning to the original (depth)
coordinate,
$$
   v_{\mathrm{II}}(Z)
      \sim D \left(1 - \frac{b_2}{5 b_1} Z \right) \,
   \mathrm{Ai}\left[b_1^{1/3} h^{-2/3}
             \left(Z + \frac{b_2Z^2}{5 b_1}\right)\right]
$$
so that
\begin{multline}
   u_{\mathrm{II}}(Z) \sim D \left(E^{-1/2}
        - \frac{1}{2} E^{-3/2} a_1 Z\right)
\\
     \left(1 - \frac{a_2 E - a_1^2}{5 E a_1} Z\right)
   \mathrm{Ai}\left[\left( \frac{a_1}{E}\right)^{1/3} h^{-2/3}
     \left(Z + \frac{a_2 E - a_1^2}{5 E a_1} Z^2 \right)\right] .
\end{multline}

Now, we examine $u_{\mathrm{I}}(Z)$ for small $Z$. We make the
following approximations
\begin{eqnarray*}
   [\hat{\mu}(Z) (\hat{\mu}(Z) - E)]^{-1/4} &\sim&
       Z^{-1/4} (E a_1)^{-1/4}
   \left( 1 - \frac{1}{4} \frac{E a_2 + a_1^2}{E a_1} Z \right) ,
\\
   \int_0^Z \sqrt{\frac{\hat{\mu} - E}{\hat{\mu}}} \mathrm{d}Z'
   &\sim& \frac{2}{3} Z^{3/2} \left(\frac{a_1}{E}\right)^{1/2}
          + \frac{E a_2 - a_1^2}{5 E a_1}
              \left( \frac{a_1}{E} \right)^{1/2} Z^{5/2} ,
\\
   \mathcal{S}_2(\delta,Z) &\sim&
   - \frac{5}{48} E^{1/2} a_1^{-1/2} Z^{-3/2}
              - \frac{E^{1/2} a_2 a_1^{-3/2}}{12} \delta^{-1/2} .
\end{eqnarray*}
In the asymptotic expansion of $\mathcal{S}_2$, we neglect terms
$\mathcal{O}(Z^{-1/2})$, which is justified because $h Z^{-1/2}$ is
small (compared to $h \delta^{-1/2}$, $h Z^{-3/2}$) in the limit $h \to
0$. Substituting these formulas into $u_{\mathrm{I}}$ gives
\begin{multline*}
   u_{\mathrm{I}} \sim C Z^{-1/4} (E a_1)^{-1/4}
      \left(1 - \frac{1}{4} \frac{E a_2 + a_1^2}{E a_1} Z\right)
\\
   \exp\left[-\frac{2}{3h} Z^{3/2} \left(\frac{a_1}{E}\right)^{1/2}
   - \frac{1}{5h} \frac{E a_2 - a_1^2}{E a_1}
     \left(\frac{a_1}{E}\right)^{1/2} Z^{5/2}\right.
\\
   \left.
   - \frac{h}{48} E^{1/2} a_1^{-5/2} Z^{-3/2}
   - \frac{h}{12} a^{-3/2} (2 a_1 - E) E^{-1/2}Z^{-3/2}
   - \frac{h E^{1/2} a_2 a_1^{-3/2}}{12} \mu^{-1/2}\right] .
\end{multline*}
Revisiting $u_{\mathrm{II}}(Z)$, for large positive $Z$, we employ the
asymptotic behavior of $\mathrm{Ai}$,
\[
   \mathrm{Ai}(s) \sim \frac{1}{2\sqrt{\pi}} s^{-1/4}
   \left(1 - \frac{5}{48} s^{-3/2}\right)
                               \exp\left[-\frac{2}{3} s^{3/2}\right]
\]
and obtain
\begin{multline*}
   u_{\mathrm{II}}(Z) \sim D \frac{1}{2 \sqrt{\pi}}
   \left(\frac{a_1}{E}\right)^{-1/12} h^{1/6} Z^{-1/4} E^{-1/2}
   \left(1 - \frac{E a_2 + a_1^2}{4 E a_1} Z\right)
\\
   \left(1 - \frac{5}{48} h Z^{-3/2} \left(\frac{a_1}{E}\right)^{-1/2}
         \right)
   \exp\left[-\frac{2}{3} \left(\frac{a_1}{E}\right)^{1/2} h^{-1}
      Z^{3/2} \left(1 + \frac{3}{2} \left(\frac{a_2 E - a_1^2}{5 E a_1}
              \right) Z \right)\right] .
\end{multline*}
Uniformly asymptotically matching $u_{\mathrm{I}}$ and
$u_{\mathrm{II}}$ then leads to the condition,
\begin{equation}
   C = \frac{D}{2\sqrt{\pi}} h^{1/6}
   \exp\left[\frac{h E^{1/2} a_2 a_1^{-3/2}}{12} \mu^{-1/2}\right]
             a_1^{1/6} E^{-1/3} .
\end{equation}

In \textbf{region III}, we construct the (oscillatory) WKB solution,
\begin{multline} \label{eq14}
   u_{\mathrm{III}}(Z) \sim F [(E - \hat{\mu}) \hat{\mu}]^{-1/4}
   \exp\left[\frac{\mathrm{i}}{h} \mathcal{S}_0(Z)
                 + \mathrm{i} h \mathcal{S}_2(\mu,Z)\right]
\\
   + G [(E - \hat{\mu}) \hat{\mu}]^{-1/4}
                 \exp\left[-\frac{\mathrm{i}}{h}
       \mathcal{S}_0(Z) - \mathrm{i} h \mathcal{S}_2(\mu,Z)\right] ,
   \quad\quad Z,h\rightarrow 0^+ ,
\end{multline}
where
\begin{equation} \label{eq:calS0b}
   \mathcal{S}_0(Z) = \int_{Z}^0
         \sqrt{\frac{E - \hat{\mu}}{\hat{\mu}}} \mathrm{d}Z'
\end{equation}
and
\begin{multline} \label{eq:calS2b}
   \mathcal{S}_2(\delta,Z) = -\int_{Z}^{-\delta}
   \left[\frac{(E \hat{\mu}' - 2 \hat{\mu} \hat{\mu}')^2}{
          32 \hat{\mu}^{3/2} (E - \hat{\mu})^{5/2}}
       + \frac{-E^2 \hat{\mu}'' + 3 E \hat{\mu} \hat{\mu}''
                    - 2 \hat{\mu}^2 \hat{\mu}''
         + E (\hat{\mu}')^2}{8 (E - \hat{\mu})^{5/2}
                    \hat{\mu}^{1/2}}\right] \mathrm{d}Z'
\\
   = \frac{(3 E + 2 \hat{\mu}) \hat{\mu}'}{
              48 \hat{\mu}^{1/2} (E - \hat{\mu})^{3/2}}
       + \frac{\hat{\mu}'}{24 (E - \hat{\mu})^{1/2} \hat{\mu}^{1/2}}
     + \int_{Z}^{-\delta}
         \left[\frac{(\hat{\mu}')^2}{
                24 \hat{\mu}^{3/2} (E - \hat{\mu})^{1/2}}
     - \frac{(7 E - 8 \hat{\mu}) \hat{\mu}''}{
                48 \hat{\mu}^{1/2} (E - \hat{\mu})^{3/2}}
                             \right] \mathrm{d}Z' .
\end{multline}
Next, we uniformly asymptotically match $u_{\mathrm{II}}$ and
$u_{\mathrm{III}}$. To this end, we consider the asymptotic behavior
of $\mathrm{Ai}(s)$ for large negative $s$,
\[
   \mathrm{Ai}(s) \sim \frac{1}{\sqrt{\pi}} s^{-1/4}
       \sin\left[-\frac{2}{3} s^{3/2} + \frac{\pi}{4}\right] ,
\]
and obtain
\begin{multline*}
   u_{\mathrm{II}}(Z) \sim D \frac{1}{\sqrt{\pi}}
     \left(\frac{a_1}{E}\right)^{-1/12} h^{1/6} Z^{-1/4} E^{-1/2}
       \left(1 - \frac{E a_2 + a_1^2}{4 E a_1} Z\right)
\\
   \sin\left[-\frac{2}{3} \left(\frac{a_1}{E}\right)^{1/2} h^{-1}
     Z^{3/2} \left(1 +
     \frac{3}{2} \left(\frac{a_2 E - a_1^2}{5 E a_1}\right) Z\right)
             + \frac{\pi}{4}\right] ,\quad\quad Z,h\rightarrow 0^+ .
\end{multline*}
Matching requires that $u_{\mathrm{III}}(Z)$ has form
\begin{multline*}
   u_{\mathrm{III}}(Z) \sim \frac{D}{\sqrt{\pi}}
          \left(\frac{a_1}{E}\right)^{-1/12} h^{1/6} E^{-1/2}
                                [(E - \hat{\mu}) \hat{\mu}]^{-1/4}
\\
   \sin\left[\frac{1}{h} \mathcal{S}_0(Z)
        + \frac{\pi}{4} + h \mathcal{S}_2(\delta,Z)
        - \frac{h E^{1/2} a_2 a_1^{-3/2}}{12} \delta^{-1/2}\right] ,
   \quad\quad Z,h\rightarrow 0^+ .
\end{multline*}
Thus
\begin{eqnarray}
   F &=&\frac{D}{2 \sqrt{\pi}}
        \left(\frac{a_1}{E}\right)^{-1/12} h^{1/6} E^{-1/2}
     \exp\left[ \frac{\mathrm{i} \pi}{4}
         - \mathrm{i} \frac{h E^{1/2} a_2 a_1^{-3/2} \delta^{-1/2}}{
                            12} \right] ,
\\
   G &=&-\frac{D}{2\sqrt{\pi}}
         \left(\frac{a_1}{E}\right)^{-1/12} h^{1/6} E^{-1/2}
     \exp\left[ -\frac{\mathrm{i} \pi}{4}
         + \mathrm{i} \frac{h E^{1/2} a_2 a_1^{-3/2} \delta^{-1/2}}{
                            12} \right] .
\end{eqnarray}

\medskip\medskip

The \textbf{Neumann boundary condition} pertains to region III, is
applied at $Z = Z_E$ in the shifted coordinate and yields the
Bohr-Sommerfeld rule. It takes the implicit form
\begin{multline} \label{eq:BS-1}
   \cot\left[ \frac{1}{h} \mathcal{S}_0(Z_E)
           + \frac{\pi}{4} + h \mathcal{S}_2(\delta,Z_E)
        - \frac{h E^{1/2} a_2 a_1^{-3/2}}{24} \delta^{-1/2} \right]
\\
   = \mathfrak{F}(h,E) ,\quad
     \mathfrak{F}(h,E) = \left.\frac{h (E - 2 \hat{\mu})
              \hat{\mu}'}{4 (E - \hat{\mu}) \hat{\mu}
     \left(-\sqrt{\frac{E - \hat{\mu}}{\hat{\mu}}}
              + h^2 \mathcal{S}_2'\right)} \right|_{Z = Z_E} .
\end{multline}
We carry out an asymptotic expansion of
$\cot^{-1} (\mathfrak{F}(h,E))$ in the small $h$ limitm
\[
   \cot^{-1}(\mathfrak{F}(h,E))
      = \frac{\pi}{2} + h \mathfrak{F}_1(E) + \mathcal{O}(h^2) ,
\]
where
\[
   \mathfrak{F}_1(E) = \left.
     \frac{(E - 2 \hat{\mu}) \hat{\mu}'}{4 (E - \hat{\mu})^{3/2}
               \hat{\mu}^{1/2}}\right|_{Z = Z_E} .
\]

We undo the shift and return to the original (depth) coordinate. We
consider, again, a function $f$ such that $\hat{\mu}(f(E)) = E$ when
$Z_E = f(E)$. Substituting (\ref{eq:calS0b}) and (\ref{eq:calS2b}),
(\ref{eq:BS-1}) takes the form
\begin{multline*}
   \frac{1}{h} \int_{f(E)}^0 \sqrt{\frac{E - \hat{\mu}}{\hat{\mu}}}
                             \mathrm{d}Z
   + \frac{\pi}{4} + \frac{(3 E + 2 \hat{\mu}(0)) \hat{\mu}'(0)}{
               48 \hat{\mu}^{1/2} (E - \hat{\mu}(0))^{3/2}}
   + \frac{\hat{\mu}'(0)}{
               24 (E - \hat{\mu}(0))^{1/2} \hat{\mu}^{1/2}(0)}
\\
   + \int_{f(E) + \delta}^{0} \left[
     \frac{(\hat{\mu}')^2}{24 \hat{\mu}^{3/2}(E - \hat{\mu})^{1/2}}
   - \frac{(7 E - 8 \hat{\mu}) \hat{\mu}''}{
               48 \hat{\mu}^{1/2} (E - \hat{\mu})^{3/2}} \right]
          \mathrm{d}Z
   - \frac{h E^{1/2} a_2 a_1^{-3/2}}{24} \delta^{-1/2}
\\
   = \left(\alpha - \frac{1}{2}\right) \pi + h \mathfrak{F}_1(E) ,
   \quad \alpha = 1, 2, \cdots,
\end{multline*}
where
\[
   \mathfrak{F}_1(E)
       = \frac{(-E + 2 \hat{\mu}(0)) \hat{\mu}'(0)}{
                  4 (E - \hat{\mu}(0))^{3/2} \hat{\mu}^{1/2}(0)} .
\]
By letting $\delta \downarrow 0$, using that
\[
   \hat{\mu}(f(E) + \delta) - E \sim -a_1 \delta + a_2 \delta^2 ,
\]
where $a_1 > 0$, and that
\[
   \frac{E^{1/2} a_2 a_1^{-3/2}}{24} \delta^{-1/2}
   \sim -\frac{E \hat{\mu}''(f(E))}{
       12 \sqrt{\hat{\mu}(f(E)) (E - \hat{\mu}(f(E) - \delta))}}
                     \frac{1}{\hat{\mu}'(E)} ,
\]
we obtain the quantization rule,
\begin{equation*}
   \frac{1}{h} \frac{1}{4}\widetilde{S}_0(E) + \frac{\pi}{4}
           + h \frac{1}{4}\widetilde{S}_2(E)
      = \left(\alpha - \frac{1}{2}\right) \pi + \mathcal{O}(h^2) ,
\end{equation*}
where
\begin{equation} \label{tSfirst}
   \widetilde{S}_0(E)
            = 4\int_{f(E)}^0 \sqrt{\frac{E - \hat{\mu}}{\hat{\mu}}}
                            \mathrm{d}Z
\end{equation}
and
\begin{equation}
   \frac{1}{4}\widetilde{S}_2(E)
   = \frac{(3 E + 2 \hat{\mu}(0)) \hat{\mu}'(0)}{
         48 \hat{\mu}^{1/2}(0) (E - \hat{\mu}(0))^{3/2}}
   + \frac{\hat{\mu}'(0)}{24 (E - \hat{\mu}(0))^{1/2}
                              \hat{\mu}^{1/2}(0)}
     - \frac{1}{24} \frac{\mathrm{d}}{\mathrm{d}E} \widetilde{J}(E)
     - \frac{1}{8} \widetilde{K}(E) - \mathfrak{F}_1(E) ,
\end{equation}
in which
\begin{eqnarray}
   \widetilde{J}(E) &=& \int_{f(E)}^0
        \left( E \hat{\mu}'' - \frac{2 (E - \hat{\mu})}{\hat{\mu}}
               (\hat{\mu}')^2 \right)
             \frac{\mathrm{d}Z}{\sqrt{\hat{\mu} (E - \hat{\mu})}} ,
\\
   \widetilde{K}(E) &=& \int_{f(E)}^0
   \hat{\mu}'' \frac{\mathrm{d}Z}{\sqrt{\hat{\mu} (E - \hat{\mu})}} .
\label{tSlast}
\end{eqnarray}
This quantization rule is satisfied by $E = \nu_{\alpha}(h)$ for the
half well.

\medskip

\begin{remark}
The above quantization rule suggest that $\lambda_1 = E_0 +
\mathcal{O}(h^{2/3})$ under Assumption \ref{assu1}, since the first eigenvalue is (semiclassically) associated with the half well. This would give us an improved version of
Lemma~\textnormal{\ref{first}}. If $\hat{\mu}'(0) = 0$, then the same
quantization rule would lead to $\lambda_1 = E_0 + \mathcal{O}(h)$.
\end{remark}

\subsection{Full well}

In anticipation that the Neumann boundary condition will not play a
role, here, we consider the eigenvalue problem on the entire real
line. We assume that there are two simple turning points, at $Z =
f_-(E)$ and at $Z = f_+(E)$; that is, $\hat{\mu} < E$ on
$]f_-(E),f_+(E)[$, and $\hat{\mu} > E$ on $]-\infty,f_-(E)[$ and
$]f_+(E),+\infty[$. Clearly, $\hat{\mu}(f_-(E)) = \hat{\mu}(f_+(E)) =
E$. Similar to the half well case, now, we construct WKB solutions in
the different regions and match them in the neighborhoods of the two
turning points $f_-(E)$ and $f_+(E)$. We let $a_{1,-}, a_{2,-}$ and
$a_{1,+}, a_{2,+}$ be the expansion coefficients of $\hat{\mu} - E$ in
the neighborhoods of $f_-(E)$ and $f_+(E)$, respectively. We now have
\begin{equation*}
\begin{split}
& \lim_{\delta \downarrow 0}
  \int_{f_-(E) + \delta}^{f_+(E) - \delta}
  -\frac{(7 E - 8 \hat{\mu}''}{
        48 \hat{\mu}^{1/2} (E - \hat{\mu})^{3/2}} \mathrm{d}Z
      - \frac{E^{1/2} a_{2,-} a_{1,-}^{-3/2}}{12} \delta^{-1/2}
      - \frac{E^{1/2} a_{2,+} a_{1,+}^{-3/2}}{12} \delta^{-1/2}
\\
=& \lim_{\delta \downarrow 0}
   \int_{f_-(E) + \delta}^{f_+(E) - \delta}
   \left( -\frac{\hat{\mu}''}{
        24 \hat{\mu}^{1/2}(E - \hat{\mu})^{1/2}}
   + \frac{E \hat{\mu}''}{48 \hat{\mu}^{1/2}(E - \hat{\mu})^{3/2}}
          \right) \mathrm{d}Z
\\
& \quad\quad
  + \frac{E \hat{\mu}''(f_-(E))}{24 \sqrt{\hat{\mu}(f_-(E))
   (E - \hat{\mu}(f_-(E) + \delta))}} \frac{1}{\hat{\mu}'(f_-(E))}
\\
& \quad\quad
  - \frac{E \hat{\mu}''(f_+(E))}{24 \sqrt{\hat{\mu}(f_+(E))
   (E - \hat{\mu}(f_+(E) - \delta))}} \frac{1}{\hat{\mu}'(f_+(E))}
\\
& \quad\quad
  - \int_{f_-(E)}^{f_+(E)}
     \frac{\hat{\mu}''}{8 \hat{\mu}^{1/2}(E - \hat{\mu})^{1/2}}
                            \mathrm{d}Z
  + \int_{f_-(E)}^{f_+(E)}
     \frac{(\hat{\mu}')^2}{
           24 \hat{\mu}^{3/2}(E - \hat{\mu})^{1/2}} \mathrm{d}Z
\\
=& -\frac{1}{24} \frac{\mathrm{d}}{\mathrm{d}E}J(E)
                    - \frac{1}{8}K(E) ,
\end{split}
\end{equation*}
where
\begin{eqnarray}
   J(E) &=& \int_{f_{-}(E)}^{f_{+}(E)}
   \left( E \hat{\mu}'' - \frac{2 (E - \hat{\mu})}{\hat{\mu}}
          (\hat{\mu}')^2 \right)
             \frac{\mathrm{d}Z}{\sqrt{\hat{\mu} (E - \hat{\mu})}} ,
\label{eq:JEsw}\\
   K(E) &=& \int_{f_{-}(E)}^{f_{+}(E)} \hat{\mu}''
             \frac{\mathrm{d}Z}{\sqrt{\hat{\mu} (E - \hat{\mu})}} .
\label{eq:KEsw}
\end{eqnarray}
That is, we arrive at the quantization
\begin{equation*}
   \frac{1}{h} \frac{1}{2}S_0(E) + h \frac{1}{2}S_2(E)
          \sim \left(\alpha - \frac{1}{2}\right) \pi ,
\end{equation*}
where
\begin{equation} \label{formSfirst}
   S_0(E) = \frac{1}{2}\int_{f_-(E)}^{f_+(E)}
            \sqrt{\frac{E - \hat{\mu}}{\hat{\mu}}} \mathrm{d}Z
\end{equation}
and
\begin{equation} \label{formSsecond}
   S_2(E) = -\frac{1}{12} \frac{\mathrm{d}}{\mathrm{d}E} J(E)
                   - \frac{1}{4} K(E) .
\end{equation}
This quantization rule is satisfied by $E = \mu_{\alpha}(h)$ for the
full well. We note that the above form has also been derived in
\cite{CdV2011} using the method introduced in \cite{CdV2005}.

\subsection{Multiple wells}

In the case of multiple wells we invoke

\medskip

\begin{assumption} \label{assu_boundary well}
There is a $Z^* < 0$ such that $\hat{\mu}'(Z^*) = 0$,
$\hat{\mu}''(Z^*) < 0$ and $\hat{\mu}'(Z) < 0$ for $Z \in
\,]Z^*,0[\,$.
\end{assumption}

\medskip

\begin{assumption} \label{assu2jag}
The function $\hat{\mu}(Z)$ has non-degenerate critical values at a
finite set
$$\{ Z_1, Z_2, \cdots, Z_M \}$$
in $]Z_I,0[$ and all critical points are non-degenerate extrema. None
of the critical values of $\hat{\mu}(Z)$ are equal, that is,
$\hat{\mu}(Z_j) \not= \hat{\mu}(Z_k)$ if $j \not= k$.
\end{assumption}

\medskip

\noindent
We label the critical values of $\hat{\mu}(Z)$ as $E_1 < \ldots < E_M
< \hat{\mu}_I$ and the corresponding critical points by
$Z_1,\cdots,Z_M$. We use the fact that $\hat{\mu}(0) = \inf_{Z \leq 0}
\hat{\mu}(Z)$ and denote $Z_0 = 0$ and $E_0 = \hat{\mu}(Z_0)$.

\begin{figure}[htbp]
\centering
\includegraphics[width=2.5 in]{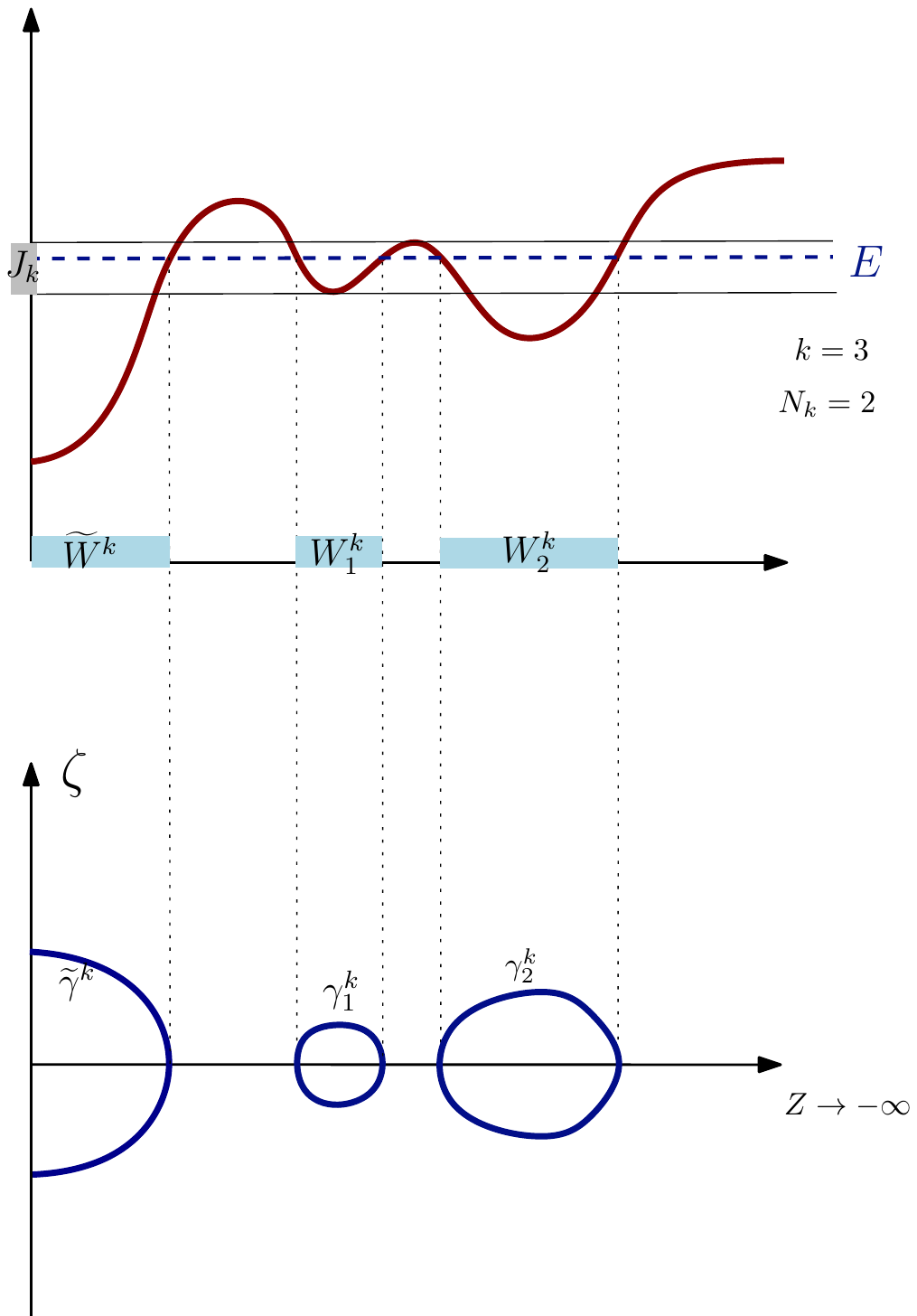}~~~~~~~~~~~~~
\includegraphics[width=2.5 in]{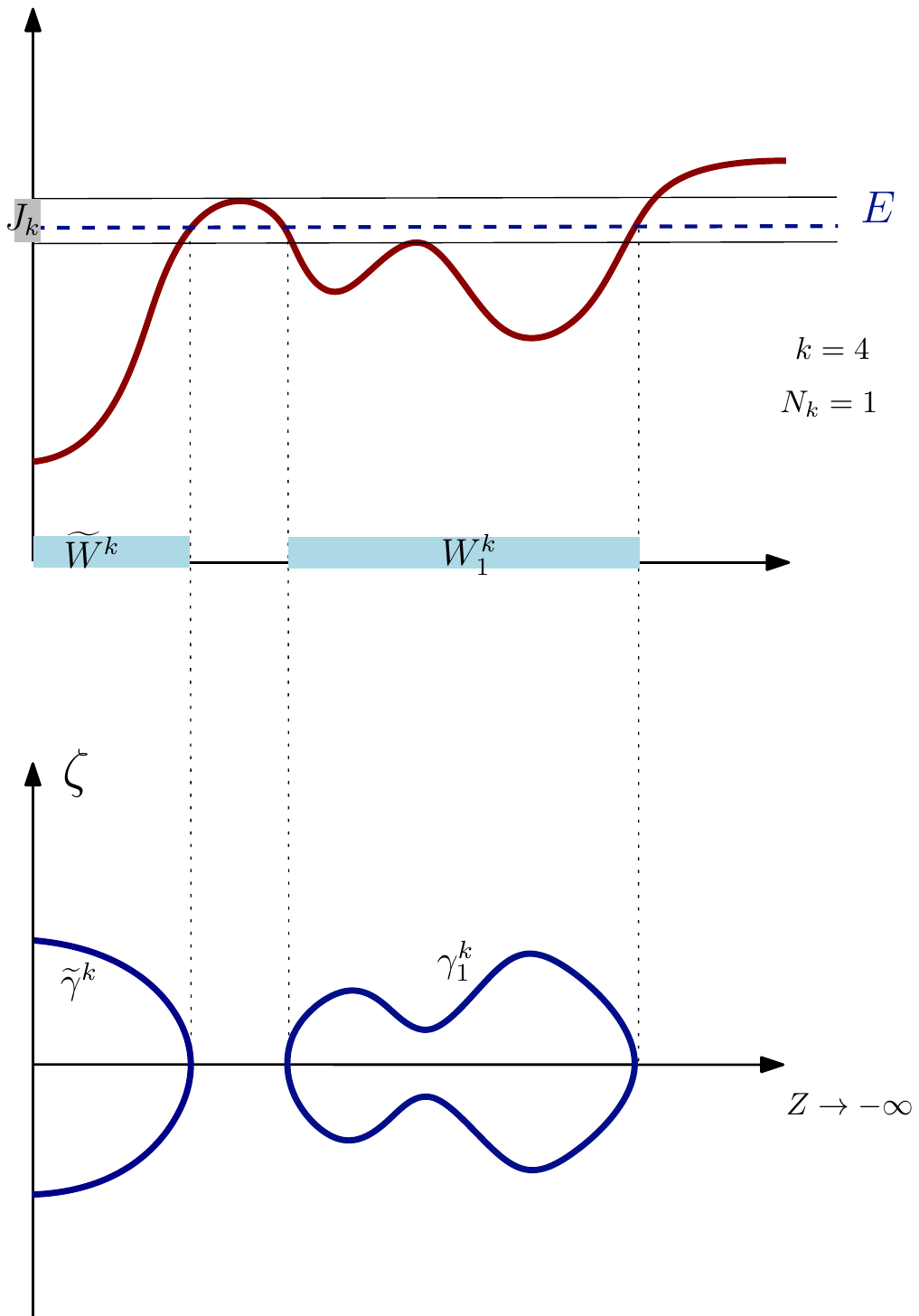}
\caption{Wells of different orders and periodic trajectories.}
\label{trajectory2}
\end{figure}

We define a well of order $k$ as a connected component of $\{ Z \in
(Z_I,0)\ :\ \hat{\mu}(Z) < E_k \}$ that is not connected to the
boundary, $Z = 0$. We refer to the connected component connected to
the boundary as a half well of order $k$. We denote $J_k =
]E_{k-1},E_k[$, $k = 1,2,3,\cdots$ and let $N_k$ ($\leq k$) be the
number of wells of order $k$ (see Figure~\ref{trajectory2} top). The
set $\{Z \in (Z_I,0)\ :\ \hat{\mu}(Z) < E_k\}$ consists of $N_k$ wells
and one half well
\[
   W_j^k(E) ,\quad j=1,2,\cdots,N_k ,\text{ and }
   \widetilde{W}^k(E) ,\quad
   (\cup_{j=1}^{N_k} W_j^k(E)) \cup \widetilde{W}^k(E)
                                       \subset [Z_I,0[ \, .
\]
The half well $\widetilde{W}^k(E)$ is connected to the boundary $Z =
0$.

Similar to Proposition 10.1 in \cite{CdV2011}, we have

\medskip

\begin{proposition}\lb{CdV Pr.10.1bis}
The semiclassical spectrum mod $o(h^{5/2})$ in $J_k$ is the union of
$N_k + 1$ spectra: $\cup_{j=1}^{N_k} \Sigma_j^k(h) \cup
\widetilde{\Sigma}^k(h)$. Here, $\Sigma_j^k(h)$ is the semiclassical
spectrum associated to well $W_j^k$, and $\widetilde{\Sigma}^k(h)$ is
the semiclassical spectrum for half well $\widetilde{W}^k$.
\end{proposition} 

\medskip

The above separation of semiclassical spectra comes from the fact that
the eigenfunctions are $\mathcal{O}(h^\infty)$ outside the wells, and
is related to the exponentially small ``tunneling'' effects
\cite{Helffer, Zelditch}. We refer further to \cite{Bender} for more
details. Therefore, we have Bohr-Sommerfeld rules for separated wells,
that is, 
\begin{equation} \label{sepa}
   \Sigma_j^k(h) = \{\mu_\alpha(h)\ :\ E_{k-1} < \mu_\alpha(h)
          < E_k\text{ and }S^{k,j}(\mu_\alpha(h)) = 2\pi h \alpha\} ,
\end{equation}
where $S^{k,j} = S^{k,j}(E) :\ ]E_{k-1},E_k[ \to \mathbb{R}$ admits
the asymptotics in $h$
$$
   S^{k,j}(E) = S^{k,j}_0(E) + h \pi + h^2 S_2^{k,j}(E)
                + \cdots 
$$
and
\begin{equation}
   \widetilde{\Sigma}^k(h) = \{\nu_\alpha(h)\ :\
   E_{k-1} < \nu_\alpha(h) < E_k\text{ and }
                    \widetilde{S}^k(\nu_\alpha(h)) = 2\pi h \alpha\},
\end{equation}
where $\widetilde{S}^k = \widetilde{S}^k(E) :\ ]E_{k-1},E_k[ \to
\mathbb{R}$ admits the asymptotics
$$
   \widetilde{S}^k(E) = \frac{1}{2} \tilde{S}^k_0(E)
        + \frac{3}{2} h \pi + \frac{1}{2} h^2 \widetilde{S}_2^k(E)
        + \cdots .
$$
The form of $S^{k,j}$ is similar to the one given in
(\ref{formSfirst})-(\ref{eq:KEsw}) and the form of $\widetilde{S}^k$
is similar to the one given in (\ref{tSfirst})-(\ref{tSlast}). We will
give more details below.

For alternative representations of $S^{k,j}$ and $\widetilde{S}^k$, we
introduce the classical Hamiltonian $p_0(Z,\zeta) = \hat{\mu}(Z) (1 +
\zeta^2)$. For any $k$, $p_0^{-1}(J_k)$ is a union of $N_k$
topological annuli $A_j^k$ and a half annulus $\widetilde{A}^k$. The
map $p_0 :\ A_j^k \to J_k$ is a fibration whose fibers $p_0^{-1}(E)
\cap A_j^k$ are topological circles $\g_j^k(E)$ that are periodic
trajectories of classical dynamics (illustrated in
Figure~\ref{trajectory2} bottom). The map $p_0 :\ \widetilde{A}^k \to
J_k$ is a topological half circle $\widetilde{\gamma}^k(E)$. If $E \in
J_k$ then $p_0^{-1}(E) = (\cup_{j=1}^{N_k} \g_j^k(E)) \cup
\widetilde{\gamma}^k(E)$. The corresponding classical periods are
\[
   T_j^k(E) = \int_{\g_j^k(E)}|\mathrm{d}t| .
\]
We let $t$ be the parametrization of $\g^k_j(E)$ by time evolution in
\begin{equation}
\label{traj}
   \frac{\mathrm{d}Z}{\mathrm{d}t} = \partial_\zeta p_0 ,\quad
   \frac{\mathrm{d}\zeta}{\mathrm{d}t} = -\partial_Z p_0
\end{equation}
for a realized energy level $E$.

\begin{figure}[htbp]
\centering
\includegraphics[width=2.5 in]{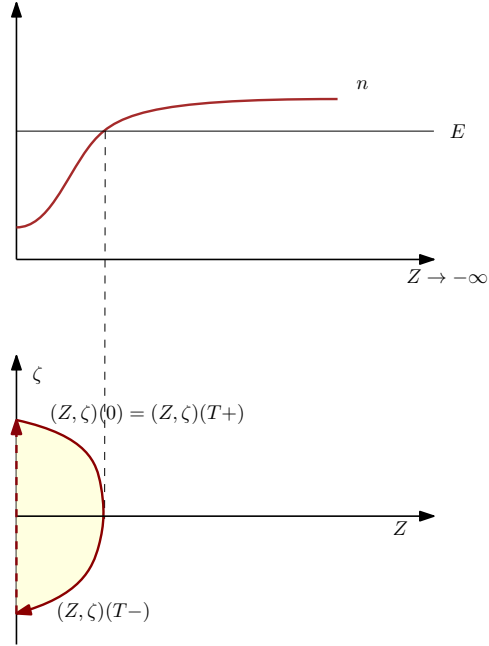}
\caption{Behavior of a half trajectory.}
\label{trajectory}
\end{figure}

For the half well $\widetilde{W}_k$, $(Z,\zeta)$ follows a periodic
(half) trajectory as shown in Figure~\ref{trajectory}. After one
(half-) period $T$, the trajectory reaches the boundary $Z(T) = 0$,
and encounters a perfect reflection, so that
\[
   \zeta(T +) = -\zeta(T -)
              = \sqrt{\frac{E - \hat{\mu}(0)}{\hat{\mu}(0)}} ,
\]
and then continues following the Hamilton system (\ref{traj}).

\subsubsection{Wells separated from the boundary}

For a well $W^k_j$ separated from the boundary, the associated
semiclassical spectrum mod $o(h^{5/2})$ follows from (\ref{sepa}) and
(\ref{eq:JEsw})-(\ref{formSsecond}). We have
\begin{equation} \label{Actions}
   S^{k,j}(E) = S_0^{k,j}(E) + h \pi + h^2 S_2^{k,j}(E) ,
\end{equation}
where
\begin{equation} \label{eq:S0}
   S_0^{k,j}(E) = \int_{\g_j^k(E)} \zeta \mathrm{d}Z
    = {\rm Area} \{(Z,\zeta)\ :\ p_0(Z,\zeta) \leq E,\,
                         Z\in W^k_j \}
\end{equation}
 and
\begin{equation} \label{eq:S2}
   S_2^{k,j}(E) = -\frac{1}{12} \frac{\mathrm{d}}{\mathrm{d}E}
        \int_{\g^k_j(E)} \left( E \hat{\mu}''
   - 2 \frac{(E - \hat{\mu})}{\hat{\mu}} (\hat{\mu}')^2 \right)
                         |\mathrm{d}t|
      - \frac14 \int_{\g^k_j(E)} \hat{\mu}'' |\mathrm{d}t| .
\end{equation}
The explicit forms of $S_0^{k,j}$ and $S_2^{k,j}$ are equivalent to
those given in (\ref{formSfirst})-(\ref{eq:KEsw}). Here, the
integration over $]f_-(E),f_+(E)[$, $E \in [E_{k-1},E_k]$, in the $Z$
coordinate has been changed into integration along the periodic
trajectory $\gamma$. One can get the same results by using the method
in \cite{CdV2005, CdV2011}. From (\ref{eq:S0}) it is immediate that
\begin{equation} \label{eq:TkjS}
   (S^{k,j}_0)'(E) = T^k_j(E) .
\end{equation}

\subsubsection{Half well connected to the boundary}
 
For the half well $\widetilde{W}^k$ connected to the boundary, we
have, mod $\mathcal{O}(h^2)$,
\begin{equation} \label{eq:tS}
   \widetilde{S}^k(E) = \frac12 \widetilde{S}^k_0(E)
             + h \frac32 \pi 
   ,
\end{equation}
where
\begin{equation} \label{eq:tS0}
   \widetilde{S}^k_0(E)
     = 2 \int_{\widetilde{\gamma}^k(E)} \zeta \mathrm{d}Z
   .
\end{equation}
The explicit form of $\widetilde{S}^k_0$ is equivalent to the one
given in (\ref{tSfirst}). Here, the integration over $]f(E),0[$, $E
\in [E_{k-1},E_k]$, in the $Z$ coordinate has been changed into
integration along the (half) periodic trajectory $\tilde{\gamma}$.
As before, it follows that
\begin{equation} \label{eq:tTkS}
   \frac12 (\widetilde{S}^k_0)'(E) = \frac{1}{2}\widetilde{T}^k(E) .
\end{equation}
The explicit form of $\widetilde{S}_2^k$ will not be needed in the
following and hence we omit it.

\medskip

\noindent
We note that $S_0^{k,j}$ and $\widetilde{S}^k_0$ depend only on
periodic trajectories.

\medskip

\begin{remark}
In the further analysis of the inverse problem, the explicit form of
$S_2^k$ is only needed for the wells separated from the boundary
(between two turning points) and there the formulas are exactly as in
\textnormal{\cite{CdV2011}} (on the whole line without boundary
conditions). Near the boundary (between a turning point and the
boundary) the function $\hat{\mu}$ is strictly decreasing and only
$S_0^k$ or the counting function for semiclassical eigenvalues suffice
to reconstruct the profile.
\end{remark}

\section{Unique recovery of $\hat{\mu}$ from the semiclassical
         spectrum}
\label{inverse}

\subsection{Trace formula}

The inverse problem is addressed with a trace formula as it reflects
the data.

\medskip

\begin{lemma}[\cite{CdV2011}, Lemma 11.1] \lb{CdV L11.1bis}
Let $S :\ J \to \R$ be a smooth function with $S' > 0$. Then we have
the following identity as Schwartz distributions in $J$, meaning that
the equality holds when applying both sides to a test function $\phi
\in C_0^\infty(J)$,
\begin{equation} \lb{CdV5bis}
   \sum_{\alpha \in \Z} \delta(E - S^{-1}(2\pi h \alpha))
        = \frac{1}{2\pi h}
   \sum_{m \in \Z} e^{\mathrm{i} m S(E) / h} S'(E) .
\end{equation}
\end{lemma}

\medskip

Substituting the action in (\ref{Actions}), (\ref{eq:tS}) and the
Bohr-Sommerfeld rules in (\ref{CdV5bis}) yields, on $J_k$ with
$\{\mu_\alpha(h)\}_\alpha = \cup_{j=1}^{N_k} \Sigma_j^k(h)$,
\begin{align*}
   \sum_{\alpha \in \Z}
   \delta(E - \mu_\alpha(h)) =& \frac{1}{2\pi h}
   \sum_{j=1}^{N_k} \sum_{m \in \Z}
        e^{\mathrm{i} m (S_0^{k,j}(E) h^{-1} + \pi + h S_2^{k,j}(E)
              + {\mathcal O}(h^2))}
   (({S^{k,j}_0})'(E) + h^2 ({S^{k,j}_2})'(E) + {\mathcal O}(h^3))
\\
   =& \frac{1}{2\pi h}
   \sum_{j=1}^{N_k} \sum_{m \in \Z}
        e^{\mathrm{i} m S_0^{k,j}(E) h^{-1}} e^{\mathrm{i} m \pi}
   ({S^{k,j}_0})'(E) (1 + \mathrm{i} m h S_2^{k,j}(E)
              + {\mathcal O}(h^2))
\end{align*}
and with $\{\nu_\alpha(h)\}_\alpha \subset \widetilde{\Sigma}^k(h)$,
\begin{align*}
   \sum_{\alpha \in \Z}
   \delta(E - \nu_\alpha(h)) =& \frac{1}{2\pi h} \sum_{m \in \Z}
   e^{\mathrm{i} m (\frac{1}{2} \widetilde{S}_0^{k}(E) h^{-1}
              + \frac{3}{2} \pi + h \frac{1}{2} \widetilde{S}_2^{k}(E)
              + {\mathcal O}(h^2))}
   \left( \frac{1}{2} (\widetilde{S}^{k}_0)'(E)
   + \frac{h^2}{2} (\widetilde{S}^k_2)'(E) + {\mathcal O}(h^3) \right)
\\
   =& \frac{1}{2\pi h} \sum_{m \in \Z}
   e^{\mathrm{i} m \frac{1}{2} \widetilde{S}_0^{k}(E) h^{-1}}
   e^{\mathrm{i} m \frac{3}{2} \pi} \frac{1}{2} (\widetilde{S}^{k}_0)'(E)
   \left( 1 + \mathrm{i} m h \frac{1}{2} \widetilde{S}_2^{k}(E)
              + {\mathcal O}(h^2) \right) .
\end{align*}
Using (\ref{eq:TkjS}) and (\ref{eq:tTkS}), and writing $\mu_\alpha$
for $\nu_\alpha$ in a unified notation, we obtain the trace formula in

\medskip

\begin{theorem} \lb{CdV Th.11.1boundary}
Let $\mu_\alpha(h)$ be the semiclassical spectrum modulo
$o(h^{5/2})$. As distributions on $J_k$, we have
\begin{align}
   \sum_{\alpha \in \Z} \delta(E - \mu_\alpha(h))
   =& \frac{1}{2\pi h} \sum_{j = 1}^{N_k}
      \sum_{m \in \Z} (-1)^m
      e^{\mathrm{i} m S_0^{k,j}(E) h^{-1}} T^k_j(E)
     (1 + \mathrm{i} m h S_2^{k,j}(E))
\nonumber\\
   &+ \frac{1}{2\pi h} \sum_{m \in \Z}
      e^{\mathrm{i} m \frac12 \widetilde{S}_0^k(E) h^{-1}}
      e^{\mathrm{i} m \frac32\pi} \widetilde{T}^k(E)
      \left( 1 + \mathrm{i} m h \frac12 \widetilde{S}^k_2(E) \right)
               + o(1) .
\label{CdV6boundary}
\end{align}
\end{theorem}

\medskip

\noindent
The direct way to obtain this trace formula is starting from
(\ref{trace_from_normal_modes}), that is,
\begin{equation*}
   \int_{\mathbb{R}^-} \widehat{\epsilon
   \partial_t \mathcal{G}_0}(Z,x,\omega,Z,\xi;\epsilon)
                 \mathrm{d}(\epsilon Z) \sim
      \frac{1}{2 h^2} \sum_{\alpha \in \Z} \delta(E - \mu_\alpha(h)) ,
\end{equation*}
upon substituting $E = h^2 \omega^2$. We then expand the parametrix
$(\ref{G0})$ in the WKB eigenfunctions $(\ref{WKB})$ from the previous
section.
    
We will use the notation
\begin{eqnarray}
   Z^k_{m,j}(E) &=& \frac{1}{2\pi h} (-1)^m
        e^{\mathrm{i} m S^{k,j}_0(E) h^{-1}} T^k_j(E)
                    (1 + \mathrm{i} m h S^{k,j}_2(E)) ,\quad
   j = 1,\cdots,N_k ,
\\
   Z^k_{m,N_k+1}(E) &=& \frac{1}{2\pi h} e^{\mathrm{i} m \frac32 \pi}
        e^{\mathrm{i} m \frac12 \widetilde{S}^k_0(E) h^{-1}}
        T^k_{N_k+1}(E)
   \left( 1 + \mathrm{i} m h \frac12 \widetilde{S}^k_2(E) \right) ,
\\
   T^k_{N_k+1}(E) :&=& \widetilde{T}^k(E)
\end{eqnarray}
for $m \in \Z$. To further unify the notation, we write
$$
   S^{k,N_k+1}_{0,2}(E) := \frac12 \widetilde{S}^k_{0,2}(E) .
$$
The micro-support of $Z^k_{m,j}$, $j=1,\cdots,N_k+1$, is given by the
Lagrangian submanifold
$$
   L^k_{m,j} = \{ (E,m T^k_j(E))\ :\ E \in J_k \}
$$
of $T^*J_k$ associated with phase function $mS^{k,j}_0(E)$.

\subsection{Separation of clusters and the weak transversality
            condition}
\label{ssec:sepclu}

We observe that the singular points of the counting function,
$\int_{p_0(Z,\zeta) \le E} |\mathrm{d}Z \mathrm{d}\zeta|$, are
precisely the critical values, $E_1, E_2, \cdots, E_M$, of $\hat{\mu}$
\cite[Lemma~11.1]{CdV2011} and, hence, are determined using the Weyl
asymptotics first. From the singularity at $E_k$ one can extract the
value of $\hat{\mu}''(Z_k)$. We then invoke

\medskip

\begin{assumption}
For any $k = 1,2,\cdots$ and any $j$ with $1 \leq j < l \leq N_k+1$,
{\em the classical periods (half-period if $j=N_k+1$) $T^k_j(E)$ and
  $T^k_l(E)$ are weakly transverse in $J_k$}, that is, there exists an
integer $N$ such that the $N$th derivative $(T^k_j - T^k_l)^{(N)}(E)$
does not vanish.
\end{assumption}

\medskip

\noindent
We introduce the sets
$$
   B = \{ E \in J_k\ :\ \exists j \neq l ,\quad
          T^k_j(E) = T^k_l(E) \} ,
$$
while suppressing $k$ in the notation. By the weak transversality
assumption, it follows that $B$ is a discrete subset of $J_k$.

We let the distributions $D_h(E) = \sum_{\alpha \in \Z} \delta(E -
\mu_\alpha(h))$ be given on intervals $J = J_k$ modulo $o(1)$ using
(\ref{CdV6boundary}). These distributions are determined mod $o(1)$ by
the semiclassical spectra mod $o(h^{5/2})$. We denote by $Z_h$ the
finite sum defined by the right-hand side of (\ref{CdV6boundary})
restricted to $m = 1$, that is,
$$
   Z_h(E) = \sum_{j=1}^{N_k+1} Z^k_{1,j}(E) .
$$
By analyzing the micro-support of $D_h$ and $Z_h$ \cite[Lemmas 12.2
  and 12.3]{CdV2011}, we find

\medskip

\begin{lemma}
Under the weak transversality assumption, the sets $B$ and the
distributions $Z_h$ mod $o(1)$ are determined by the distributions
$D_h$ mod $o(1)$.
\end{lemma}

\medskip

\begin{lemma}
Assuming that the $S^j$'s are smooth and the $a_j$'s do not vanish,
there is a unique splitting of $Z_h$ as a sum
$$
   Z_h(E) = \frac{1}{2\pi h} \sum_{j=1}^{N_k+1}
   (a_j(E) + h b_j(E)) e^{\mathrm{i} S^j(E)/h} + o(1) .
$$
\end{lemma}

\medskip

\noindent
It follows that the spectrum in $J_k$ mod $o(h^{5/2})$ determines the
actions $S^{k,j}_0(E)$, $S^{k,j}_2(E)$ and
$\widetilde{S}^k_0(E)$. This provides the separation of the data for
the $N_k$ wells and the half well.

\medskip

For the reconstruction of $\hat{\mu}$ from these actions, we need one
more assumption

\medskip

\begin{assumption} \label{defect}
The function $\hat{\mu}$ has a generic symmetry defect: If there exist
$X_\pm$ satisfying $\hat{\mu}(X_-) = \hat{\mu}(X_+) < E$, and for all
$N \in \mn$, $\hat{\mu}^{(N)}(X_-) = (-1)^N \hat{\mu}^{(N)}(X_+)$,
then $\hat{\mu}$ is globally even with respect to $\ha (X_+ + X_-)$ in
the interval $\{ Z\ :\ \hat{\mu}(Z) < E \}$.
\end{assumption}

\medskip

\noindent
We will carry out the reconstruction of $\hat{\mu}$ successively in
intervals $J_k$, $k = 1,\cdots,M$ and then on the interval
$[E_M,E_{M+1}]$ with $E_{M+1} = \hat{\mu}_I$.

\subsection{Reconstruction of a single well, with barrier and
            descreasing profile}

We discuss in detail the case of one local minimum for $Z < 0$ and
global minimum at $Z = 0$ ($\hat{\mu}(0) < \hat{\mu}(Z)$ $\forall$ $Z
< 0$, $\hat{\mu}'(0) \leq 0$). This means that the global minimum
occurs at $Z = 0$ and $E_1 = \hat{\mu}(Z_1)$ is the local
minimum. Then $E_2 = \hat{\mu}(Z_2)$ is attained at $Z_2 \in (Z_1,0)$
and $E_3 = \hat{\mu}_I$.

\begin{figure}[h]
\centering
\includegraphics[width=3 in]{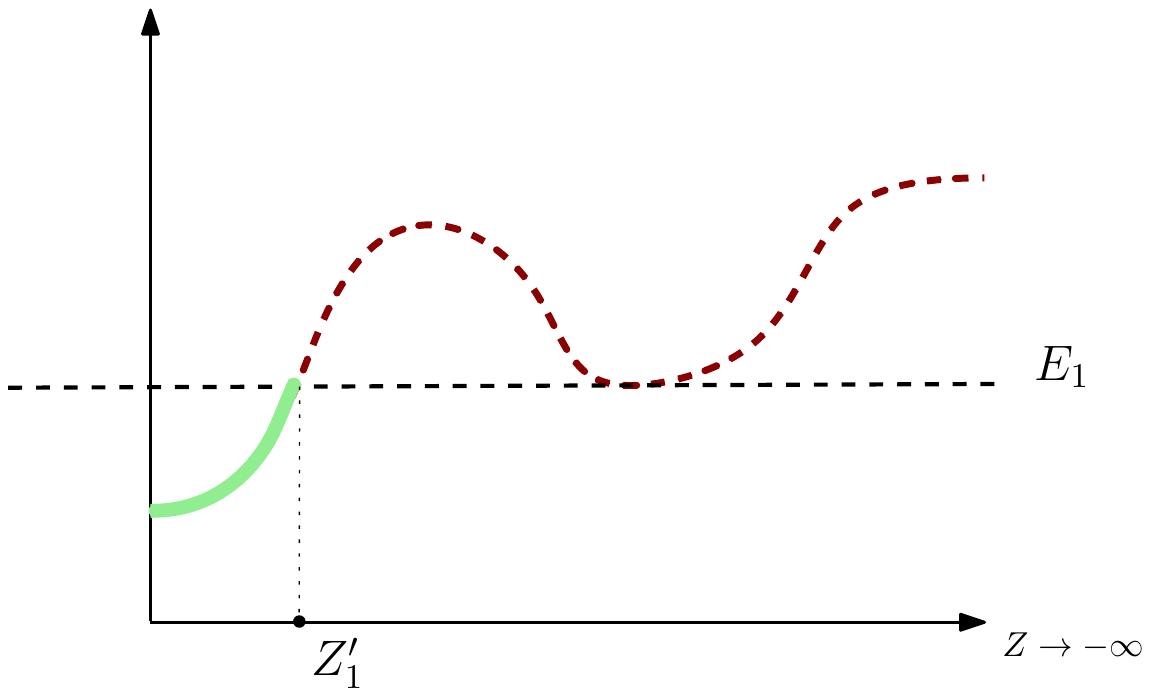}
\caption{Reconstruction Step 1 in green.}
\label{separation2}
\end{figure}

\textit{Step 1}. For $E \in ]E_0,E_1[$, there is only one (half) well,
$\widetilde{W}^1(E)$, of order $1$ with $\widetilde{W}^1(E_1) =
[Z_1',0]$. Since $\hat{\mu}$ is strictly decreasing in
$\widetilde{W}^1(E_1)$, we may reconstruct $\hat{\mu}$ on this
interval as in Section~\ref{S-decreasing}. This is illustrated in
Figure~\ref{separation2} in green.
 
\textit{Step 2}. We note that $Z_2$ in this case is the $Z^*$ defined
above Assumption~\ref{assu_boundary well}. We consider $E \in
\,]E_1,E_2[$ which corresponds to wells of order $k = 2$ with $N_k =
1$ (one connected component for $Z < 0$ separated from the
boundary). The two wells are $W^{2,1}(E)$ and $\widetilde{W}^2(E)$
with $W^{2,1}(E_2) = [Z_-,Z_2]$ and $\widetilde{W}^2(E_2) =
[Z_2,0]$. Here, $Z_-$ is the unique point in $[Z_I,Z_1]$ such that
$E_2 = \hat{\mu}(Z_-)$. We are given $S^{2,1}_0$, $S^{2,1}_2$ and
$\widetilde{S}^2_0$ (and $\widetilde{S}^2_2$).

We continue to reconstruct $\hat{\mu}$ from $[Z'_1,0]$ to $[Z_2,0]$
from $\widetilde{S}_0^2$. For the reconstruction of $\hat{\mu}$ on the
interval $I = [Z_-,Z_2]$, more effort is needed. We note that, up to
this point, $I$ itself cannot be determined yet. The following theorem
is a version of \cite[Theorem 5.1]{CdV2011}

\medskip

\begin{theorem} \label{CdV-well}
Under Assumption \ref{defect}, the function $\hat{\mu}$ is determined
on $I$ by $S_0^{2,1}$ and $S_2^{2,1}$ up to a symmetry $\hat{\mu}(Z)
\rightarrow \hat{\mu}(c - Z)$, where $\frac{c}{2}$ is the midpoint of
$I$.
\end{theorem}

\medskip

\begin{proof}
For any $E \in [E_1,E_2[$ the functions $f_\pm :\ [E_1, E_2[ \to I$,
are defined so that $W^2_1(E) = [f_-(E),f_+(E)]$. We have
$\hat{\mu}'(Z) < 0$ for $Z \in\, ]f_-(E),Z_1[$ and $\hat{\mu}'(Z) > 0$
        for $Z \in\, ]Z_1,f_+(E)[$. We introduce
\[
   \Phi(E) = f_+'(E) - f_-'(E)\quad\text{and}\quad
   \Psi(E) = \frac{1}{f_+'(E)} - \frac{1}{f_-'(E)} .
\]
As in the proof of Theorem~\ref{CdV-decr}, we have
\[
   (S^{2,1}_0)'(E) = T g(E) ,\quad
   T g(E) = \int_{E_1}^E \frac{g(u)}{\sqrt{E - u}}
                       \, \mathrm{d}u\quad\text{with}\quad
   g(u) = \frac{\Phi(u)}{\sqrt{u}} .
\]
The inversion formula for the Abel transform yields $\Phi(E)$ for $E
\in [E_1,E_2[$.

Concerning the recovery of $\Psi$, we have 
\[
   S^{2,1}_2(E) = -\frac{1}{12}\frac{\mathrm{d}}{\mathrm{d}E}\mathcal{B} \Psi(E) ,\quad
   \mathcal{B} \Psi(E)  = \int_{E_1}^E
       \left((7 E - 6 u) \Psi'(u)
     - 2 \left(\frac{E}{u} - 1 \right) \Psi(u)\right)
              \frac{\mathrm{d}u}{\sqrt{u (E - u)}} ,
\]
which follows from (\ref{formSsecond}) with
(\ref{eq:JEsw})-(\ref{eq:KEsw}) upon changing variable of integration,
$Z = f_\pm(u)$. Thus, from $S^{2,1}_2(E)$ and the fact
$\mathcal{B} \Psi(E_1) = \pi \sqrt{2
    E_1 \hat{\mu}''(Z_1)}$, we can recover $\mathcal{B} \Psi(E)$.  It
can be shown that
\[
   \frac{\pi}{E^{3/2}}
   \frac{\mathrm{d}^2}{\mathrm{d} E^2} (T \circ \mathcal{B} \Psi)(E)
   = E^2 \Psi''(E) + 4 E \Psi'(E) - \Psi(E) .
\]
That is, we obtain a second-order inhomogeneous ordinary differential
equation for $\Psi$ on the interval $[E_1,E_2[$. This equation is
    supplemented with the ``initial'' conditions 
    \[
   \Psi(E_1) = 0 ,\quad
   \lim_{E \downarrow E_1} \sqrt{E - E_1} \Psi'(E)
                            = \sqrt{2 \hat{\mu}''(Z_1)}
\]
As mentioned in Subsection~\ref{ssec:sepclu}, this second derivative
is obtained from the limiting behavior of the counting function which
coincides with $S^{2,1}_0(E)$ as $E \downarrow
E_1$. We use that the period of small
  oscillations of the ``pendulum'' associated to the local minimum of
  $\hat{\mu}$ at $Z_1$ is given by
\[
   (S^{2,1}_0)'(E) = \int_{f_-(E)}^{f_+(E)}
      \frac{\mathrm{d}Z}{\sqrt{\hat{\mu}(E - \hat{\mu})}}
   = \pi \sqrt{\frac{2}{E_1 \hat{\mu}''(Z_1)}} + o(1)\quad
  \text{as}\ E \downarrow E_1 .
\]
Thus we obtain $\Psi(E)$ for $E \in [E_1,E_2[$.

With $\pm f_\pm'(E) > 0$ for $E \in\, ]E_1,E_2[$, we then find
\begin{equation} \label{eq:fpm}
   2 f_{\pm}' = \pm \Phi + \sigma \sqrt{\Phi^2 - 4 \frac{\Phi}{\Psi}}
\end{equation}
with
\begin{equation*}
   \sigma = \sign(f_+' + f_-')
   = \left\{\begin{array}{ccc} +1 &\text{if}& f_+' + f_-' > 0
     \\
     0 &\text{if}& f_+' + f_-' = 0
     \\
     -1 &\text{if}& f_+' + f_-' < 0 \end{array}\right.
\end{equation*}
We note that the sign is not (yet) determined, and only if the well is
mirror-symmetric with respect to its vertex then $f_+' + f_-' = 0$ and
the square root in (\ref{eq:fpm}) vanishes. However,
    later, we will find the sign by a gluing argument.

By Assumption~\ref{defect}, the function $\sigma = \sigma(E)$ is
constant for $E \in\, ]E_1,E_2[$. Hence, in what follows we will
    exchange $\sigma$ with $\pm$. We have
\begin{eqnarray*}
   f_+(E) &=& Z_1 + \frac12 \int_{E_1}^E \left(
     \Phi \pm \sqrt{\Phi^2 - 4 \frac{\Phi}{\Psi}}\right) \mathrm{d}E ,
\\
   f_-(E) &=& Z_1 + \frac12 \int_{E_1}^E \left(
    -\Phi \pm \sqrt{\Phi^2 - 4 \frac{\Phi}{\Psi}}\right) \mathrm{d}E .
\end{eqnarray*}
Since $f_+(E_2) = Z_2$ and $f_-(E_2) = Z_-$, we find that
\begin{eqnarray*}
   Z_2 &=& Z_1 + \frac12 \int_{E_1}^{E_2} \left(
     \Phi \pm \sqrt{\Phi^2 - 4 \frac{\Phi}{\Psi}}\right) \mathrm{d}E ,
\\
   Z_- &=& Z_1 + \frac12 \int_{E_1}^{E_2} \left(
    -\Phi \pm \sqrt{\Phi^2 - 4 \frac{\Phi}{\Psi}}\right) \mathrm{d}E .
\end{eqnarray*}
Hence, the distance, $Z_2 - Z_1$, between the two critical points is
recovered (modulo mirror symmetry of $Z_1$ with respect to
$\frac{c}{2}$). Since $f_\pm$ are both monotonic on $]E_1,E_2[$,
$\hat{\mu}$ can be recovered (up to mirror symmetry) on $I$.
\end{proof}

\begin{figure}[h]
\centering
\includegraphics[width=3 in]{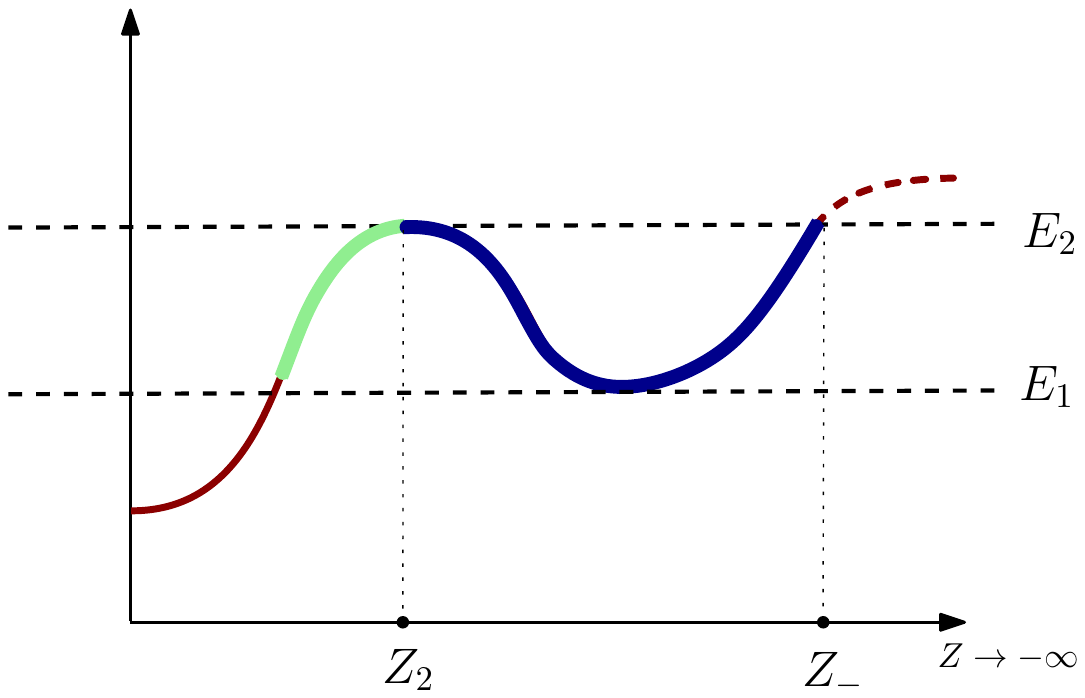}
\caption{Reconstruction Step 2, first part in green and second part in
  blue.}
\label{separation3}
\end{figure}

With this result, the reconstructions on $[Z'_1,0]$ and $I$ can be
smoothly glued together, and the uncertainty in the translation of $I$
and the ``orientation'' of $\hat{\mu}$ on $I$ are eliminated. Thus
$\hat{\mu}$ is uniquely determined on the interval $[Z_-,0]$. This is
illustrated in Figure~\ref{separation3}.

\begin{figure}[h]
\centering
\includegraphics[width=3 in]{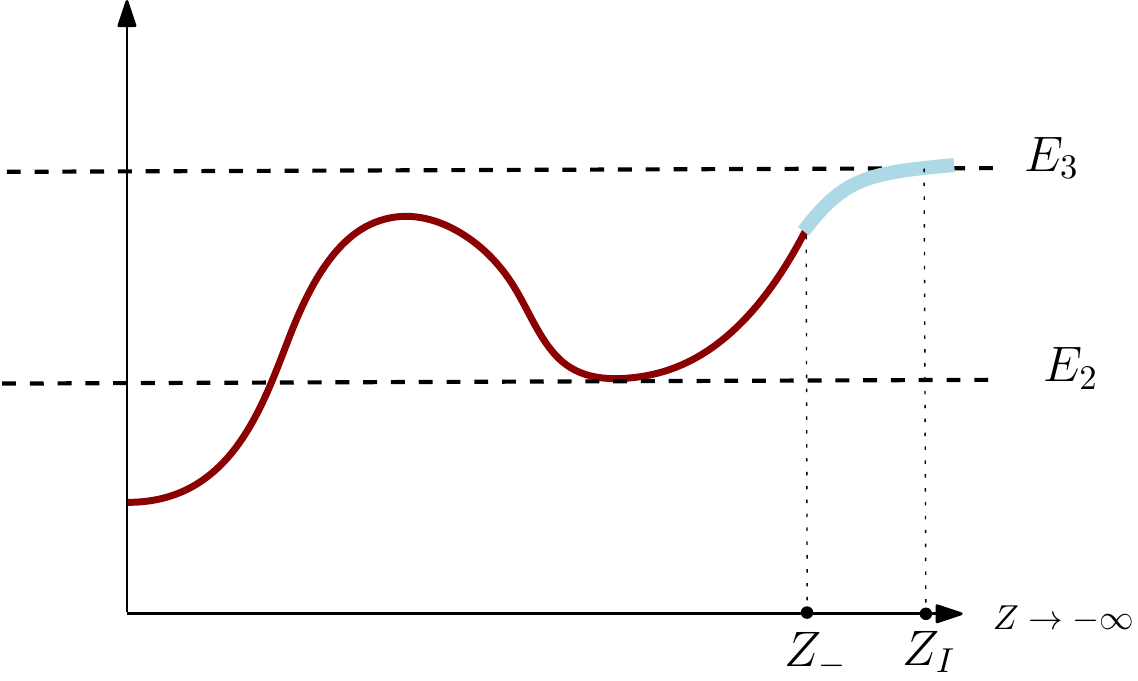}
\caption{Reconstruction Step 3 in light blue.}
\label{separation4}
\end{figure}

\textit{Step 3}. On the interval $[Z_I,Z_-]$ we may use the Weyl
asymptotics again to recover $\hat{\mu}$. The counting function in the
interval $[E_2,E_3]$ is obtained from $\widetilde{S}^3$ which
corresponds with
$$
   {\rm Area}(\{ (Z,\zeta)\ :\ \hat{\mu}(Z) (1 + \zeta^2) \leq E \})
                      = A_1(E) + A_2(E) ,
$$
where
$$
   A_1(E) = {\rm Area}(\{ (Z,\zeta)\ :\
             \hat{\mu}(Z) (1 + \zeta^2) \leq E ,\
                                     Z_- \leq Z \leq 0 \})
$$
is already known, and
$$
   A_2(E) = 2 \int_{f(E)}^{Z_-} \sqrt{\frac{E - \hat{\mu}}{\hat{\mu}}}
                                \mathrm{d}Z ,
$$
$Z_I \leq f(E) < Z_-$ since $E_2 \leq E \leq E_3 = \hat{\mu}_I$. Thus
we may recover $\hat{\mu}$ on the interval $[Z_I,Z_-]$ where
$\hat{\mu}$ is decreasing while applying Theorem~\ref{CdV-decr}. Step
3 is illustrated in Figure~\ref{separation4}.

The two profiles for $\hat{\mu}$ on $[Z_I, Z_-]$ and on $[Z_-,Z_2]$
are then glued together at $Z=Z_-$ which is already known. This
completes the reconstruction procedure.

\subsection{Reconstruction of multiple wells} 

If $\hat{\mu}$ has multiple wells, we follow an inductive
procedure. First, we consider the reconstruction of the half well
$\widetilde{W}^k$ of order $k$ between $E_{k-1}$ and $E_k$. We note
that $\widetilde{W}^k$ must be a continuation of the half well
$\widetilde{W}^{k-1}$, or be joined with some well $W^{k-1}_{j'}$ of
order $k-1$. This can be done in a fashion similar to the process
presented above (on $[Z_I,Z_-]$).

\begin{figure}[h]
\centering
\includegraphics[width=5 in]{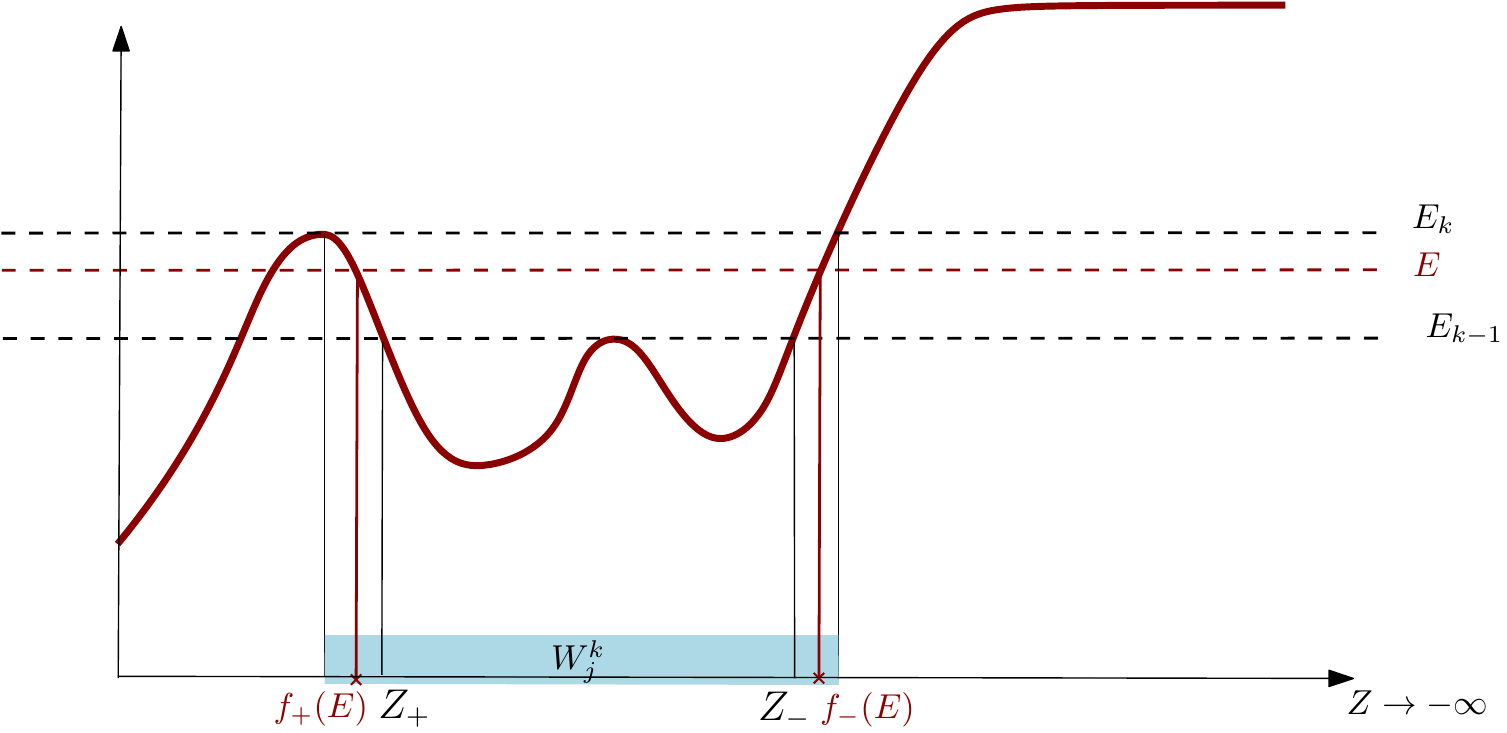}
\caption{Illustration of $f_\pm$.}
\label{fpm}
\end{figure}

Secondly, we consider the reconstruction of a well, $W^k_j$, separated
from the boundary, of order $k$. The well $W^k_j$ might be a new well,
and can be reconstructed as in Theorem \ref{CdV-well}.
The well, $W^k_j$, might also be joining two wells of order $k-1$, or
extending a single well of order $k-1$. Let the profile under
$E_{k-1}$ already be recovered. The smooth joining of two wells
can be carried out under Assumption~\ref{defect}. We consider now
functions $f_-(E)$ and $f_+(E)$ for $E \in [E_{k-1},E_{k}]$ such that
$W^k_j$ is the union of three connected intervals,
$$
   W^k_j(E_k) = [f_-(E_k),f_-(E_{k-1})[ \, \cup \,
   [f_-(E_{k-1}),f_+(E_{k-1})] \, \cup \, ]f_+(E_{k-1}),f_+(E_{k})] ;
$$
see Figure~\ref{fpm}. The semiclassical spectrum in $]E_{k-1},E_{k}[$
    up to $o(h^{5/2})$ gives the actions $S^{k,j}_0$ and $S^{k,j}_2$.

From $S^{k,j}_0$ we obtain
\begin{equation*}
   T^k_j(E) = (S^{k,j}_0)'(E)
   = \int_{f_{-}(E)}^{f_{+}(E)}
             \frac{\mathrm{d}Z}{\sqrt{\hat{\mu} (E - \hat{\mu})}}
\end{equation*}
which signify the periods of the trajectories of energy $E$. We write
$Z_- = f_-(E_{k-1})$ and $Z_+ = f_+(E_{k-1})$, and decompose the
interval,
$$
   [f_-(E),f_+(E)] = [f_-(E),Z_-] \cup [Z_-,Z_+] \cup [Z_+,f_+(E)] .
$$
In accordance with this decomposition,
\begin{equation*}
   T^k_j(E) = T_-(E) + T_{k-1}(E) + T_+(E) ,
\end{equation*}
where
\begin{eqnarray*}
   T_-(E) &=& \int_{f_{-}(E)}^{Z_-}
           \frac{\mathrm{d}Z}{\sqrt{\hat{\mu} (E - \hat{\mu})}} ,
\\
   T_{k-1}(E) &=& \int_{Z_-}^{Z_+}
           \frac{\mathrm{d}Z}{\sqrt{\hat{\mu} (E - \hat{\mu})}} ,
\\
   T_+(E) &=& \int_{Z_+}^{f_{+}(E)}
           \frac{\mathrm{d}Z}{\sqrt{\hat{\mu} (E - \hat{\mu})}} .
\end{eqnarray*}
We note that $T_{k-1}(E)$ is already known. In $T_{\mp}(E)$ we change
the variable of integration, $Z = f_{\mp}(u)$. Using that
$\hat{\mu}(f_{\mp}(u)) = u$, we get
\begin{equation*}
   T_{\mp}(E) = \mp \int_{E_{k-1}}^E
                \frac{f_{\mp}'(u)}{\sqrt{u (E - u)}} \mathrm{d}u ;
\end{equation*}
then
\begin{equation*}
   T^k_j(E) - T_{k-1}(E) = T g(E) ,\quad
   T g(E) = \int_{E_{k-1}}^E
                \frac{g(u)}{\sqrt{E - u}} \mathrm{d}u\quad
   \text{with}\quad g(u) = \frac{\Phi(u)}{\sqrt{u}}
\end{equation*}
and $\Phi(u) = f_+'(u) - f_-'(u)$ as before. Inverting this Abel
transform, we obtain $\Phi$ on $[E_{k-1},E_k[$.

From $S^{k,j}_2$ we obtain
\[
   -\frac{1}{12} \frac{\mathrm{d}}{\mathrm{d}E} J(E)
                 - \frac{1}{4} K(E) ,
\]
where
\begin{eqnarray*}
   J(E) &=& \int_{f_{-}(E)}^{f_{+}(E)}
   \left( E \hat{\mu}'' - 2 \frac{(E - \hat{\mu})}{\hat{\mu}}
         (\hat{\mu}')^2
   \right) \frac{\mathrm{d}Z}{\sqrt{\hat{\mu} (E - \hat{\mu})}} ,
\\
   K(E) &=& \int_{f_{-}(E)}^{f_{+}(E)}
   \hat{\mu}'' \frac{\mathrm{d}Z}{\sqrt{\hat{\mu} (E - \hat{\mu})}} .
\end{eqnarray*}
Using that
$$
   \hat{\mu}(f_\pm(E)) = E ,\quad
   \hat{\mu}'|_{Z = f_\pm(E)} = \frac{1}{f_\pm'(E)} ,\quad
   \hat{\mu}''|_{Z = f_\pm(E)} = \left(\frac{1}{f_\pm'}\right)'(E)
                                    \frac{1}{f_\pm'(E)} ,
$$
changing variables of integration in $J$ and $K$, $Z =
    f_\pm(u)$ and introducing
\[
   \Psi(E) = \frac{1}{f_{+}'(E)} - \frac{1}{f_{-}'(E)} ,
\]
we have
\begin{eqnarray*}
   J(E) - J_{k-1}(E) &=& \int_{E_{k-1}}^E \left(E \Psi'(u)
        - 2 \left(\frac{E}{u} - 1\right) \Psi(u)\right)
                     \frac{\mathrm{d}u}{\sqrt{u (E - u)}} ,
\\
   K(E) - K_{k-1}(E) &=& \int_{E_{k-1}}^E \Psi'(u)
                  \frac{\mathrm{d}u}{\sqrt{u (E - u)}} ,
\end{eqnarray*}
where
\begin{eqnarray*}
   J_{k-1}(E) &=& \int_{Z_-}^{Z_+} \left( E \hat{\mu}''
    - 2 \frac{(E - \hat{\mu})}{\hat{\mu}} (\hat{\mu}')^2 \right)
        \frac{\mathrm{d}Z}{\sqrt{\hat{\mu} (E - \hat{\mu})}} ,
\\
   K_{k-1}(E) &=& \int_{Z_-}^{Z_+} \hat{\mu}''
        \frac{\mathrm{d}Z}{\sqrt{\hat{\mu} (E - \hat{\mu})}}
\end{eqnarray*}
are already known. Thus, from $S^{k,j}_2$, we recover
\[
   \mathcal{B} \Psi(E) = \int_{E_{k-1}}^E
   \left( (7 E - 6 u) \Psi'(u)
      - 2 \left(\frac{E}{u} - 1\right) \Psi(u) \right)
                     \frac{\mathrm{d}u}{\sqrt{u (E - u)}} .
\]
Then similar to the proof of Theorem \ref{CdV-well}, we recover $\Psi$
on $[E_{k-1},E_k[$ by inverting $\mathcal{B}$ through the introduction
    of a second-order ordinary differential equation.

From $\Phi$ and $\Psi$ we obtain
$$
   2 f_+' = \Phi \pm \sqrt{\Phi^2 - 4 \frac{\Phi}{\Psi}} ,\quad
   2 f_-' =-\Phi \pm \sqrt{\Phi^2 - 4 \frac{\Phi}{\Psi}}
$$
and then
\begin{eqnarray*}
   f_+(E) &=& Z_+ + \frac12 \int_{E_{k-1}}^E \left(
    \Phi \pm \sqrt{\Phi^2 - 4 \frac{\Phi}{\Psi}}\right) \mathrm{d}E ,
\\
   f_-(E) &=& Z_- + \frac12 \int_{E_{k-1}}^E \left(
    \Phi \pm \sqrt{\Phi^2 - 4 \frac{\Phi}{\Psi}}\right) \mathrm{d}E .
\end{eqnarray*}
From $f_-$ we recover $\hat{\mu}$ on the interval
    $[f_-(E),Z_-]$ and from $f_+$ we recover $\hat{\mu}$ on the
    interval $[Z_+,f_+(E)]$. The $\pm$ signs in $f_{\pm}$ are
    disentangled by smoothly joining the newly reconstructed pieces to
    the previously reconstructed part and Assumption~\ref{defect},
as in previous section. Since the profile in $[Z_-,Z_+]$ can only be
determined up to translation and symmetry, the determination of the
profile in $W^k_j$ is up to the same translation and symmetry.

The symmetry and translation freedom for all the wells will be
gradually eliminated during the whole process. At the final step,
there is a single half well connected to the boundary, and then we can
reconstruct exactly the entire profile.

\section*{Acknowledgement}

The authors thank Y. Colin de Verdi\`ere for invaluable
discussions.  M.V.d.H. gratefully acknowledges support from the
Simons Foundation under the MATH $+$ X program, from NSF under grant
DMS-1559587 and the Geo-Mathematical Imaging Group at Rice University.


\end{document}